\documentclass[twoside]{report}
\usepackage[a4paper]{geometry}
\usepackage[latin1]{inputenc} 
\usepackage[T1]{fontenc} 
\usepackage{hyperref}
\usepackage{dsfont}
\usepackage{graphicx, amsmath, theorem, amssymb, algorithmic}
\usepackage{color}

\usepackage{caption}
\usepackage{subcaption}
\usepackage{tikz}
\usepackage{pgfplots}
\usetikzlibrary{spy,calc}
\usepackage{hyperref}

\newif\ifblackandwhitecycle
\gdef\patternnumber{0}

\pgfkeys{/tikz/.cd,
    zoombox paths/.style={
        draw=orange,
        very thick
    },
    black and white/.is choice,
    black and white/.default=static,
    black and white/static/.style={ 
        draw=white,   
        zoombox paths/.append style={
            draw=white,
            postaction={
                draw=black,
                loosely dashed
            }
        }
    },
    black and white/static/.code={
        \gdef\patternnumber{1}
    },
    black and white/cycle/.code={
        \blackandwhitecycletrue
        \gdef\patternnumber{1}
    },
    black and white pattern/.is choice,
    black and white pattern/0/.style={},
    black and white pattern/1/.style={    
            draw=white,
            postaction={
                draw=black,
                dash pattern=on 2pt off 2pt
            }
    },
    black and white pattern/2/.style={    
            draw=white,
            postaction={
                draw=black,
                dash pattern=on 4pt off 4pt
            }
    },
    black and white pattern/3/.style={    
            draw=white,
            postaction={
                draw=black,
                dash pattern=on 4pt off 4pt on 1pt off 4pt
            }
    },
    black and white pattern/4/.style={    
            draw=white,
            postaction={
                draw=black,
                dash pattern=on 4pt off 2pt on 2 pt off 2pt on 2 pt off 2pt
            }
    },
    zoomboxarray inner gap/.initial=5pt,
    zoomboxarray columns/.initial=2,
    zoomboxarray rows/.initial=2,
    subfigurename/.initial={},
    figurename/.initial={zoombox},
    zoomboxarray/.style={
        execute at begin picture={
            \begin{scope}[
                spy using outlines={%
                    zoombox paths,
                    width=\imagewidth / \pgfkeysvalueof{/tikz/zoomboxarray columns} - (\pgfkeysvalueof{/tikz/zoomboxarray columns} - 1) / \pgfkeysvalueof{/tikz/zoomboxarray columns} * \pgfkeysvalueof{/tikz/zoomboxarray inner gap} -\pgflinewidth,
                    height=\imageheight / \pgfkeysvalueof{/tikz/zoomboxarray rows} - (\pgfkeysvalueof{/tikz/zoomboxarray rows} - 1) / \pgfkeysvalueof{/tikz/zoomboxarray rows} * \pgfkeysvalueof{/tikz/zoomboxarray inner gap}-\pgflinewidth,
                    magnification=8,
                    every spy on node/.style={
                        zoombox paths
                    },
                    every spy in node/.style={
                        zoombox paths
                    }
                }
            ]
        },
        execute at end picture={
            \end{scope}
            \node at (image.north) [anchor=north,inner sep=0pt] {\subcaptionbox{\label{\pgfkeysvalueof{/tikz/figurename}-image}}{\phantomimage}};
            \node at (zoomboxes container.north) [anchor=north,inner sep=0pt] {\subcaptionbox{\label{\pgfkeysvalueof{/tikz/figurename}-zoom}}{\phantomimage}};
     \gdef\patternnumber{0}
        },
        spymargin/.initial=0.5em,
        zoomboxes xshift/.initial=1,
        zoomboxes right/.code=\pgfkeys{/tikz/zoomboxes xshift=1},
        zoomboxes left/.code=\pgfkeys{/tikz/zoomboxes xshift=-1},
        zoomboxes yshift/.initial=0,
        zoomboxes above/.code={
            \pgfkeys{/tikz/zoomboxes yshift=1},
            \pgfkeys{/tikz/zoomboxes xshift=0}
        },
        zoomboxes below/.code={
            \pgfkeys{/tikz/zoomboxes yshift=-1},
            \pgfkeys{/tikz/zoomboxes xshift=0}
        },
        caption margin/.initial=4ex,
    },
    adjust caption spacing/.code={},
    image container/.style={
        inner sep=0pt,
        at=(image.north),
        anchor=north,
        adjust caption spacing
    },
    zoomboxes container/.style={
        inner sep=0pt,
        at=(image.north),
        anchor=north,
        name=zoomboxes container,
        xshift=\pgfkeysvalueof{/tikz/zoomboxes xshift}*(\imagewidth+\pgfkeysvalueof{/tikz/spymargin}),
        yshift=\pgfkeysvalueof{/tikz/zoomboxes yshift}*(\imageheight+\pgfkeysvalueof{/tikz/spymargin}+\pgfkeysvalueof{/tikz/caption margin}),
        adjust caption spacing
    },
    calculate dimensions/.code={
        \pgfpointdiff{\pgfpointanchor{image}{south west} }{\pgfpointanchor{image}{north east} }
        \pgfgetlastxy{\imagewidth}{\imageheight}
        \global\let\imagewidth=\imagewidth
        \global\let\imageheight=\imageheight
        \gdef\columncount{1}
        \gdef\rowcount{1}
        
    },
    image node/.style={
        inner sep=0pt,
        name=image,
        anchor=south west,
        append after command={
            [calculate dimensions]
            node [image container,subfigurename=\pgfkeysvalueof{/tikz/figurename}-image] {\phantomimage}
            node [zoomboxes container,subfigurename=\pgfkeysvalueof{/tikz/figurename}-zoom] {\phantomimage}
        }
    },
    color code/.style={
        zoombox paths/.append style={draw=#1}
    },
    connect zoomboxes/.style={
    spy connection path={\draw[draw=none,zoombox paths] (tikzspyonnode) -- (tikzspyinnode);}
    },
    help grid code/.code={
        \begin{scope}[
                x={(image.south east)},
                y={(image.north west)},
                font=\footnotesize,
                help lines,
                overlay
            ]
            \foreach \x in {0,1,...,9} { 
                \draw(\x/10,0) -- (\x/10,1);
                \node [anchor=north] at (\x/10,0) {0.\x};
            }
            \foreach \y in {0,1,...,9} {
                \draw(0,\y/10) -- (1,\y/10);                        \node [anchor=east] at (0,\y/10) {0.\y};
            }
        \end{scope}    
    },
    help grid/.style={
        append after command={
            [help grid code]
        }
    },
}

\newcommand\phantomimage{%
    \phantom{%
        \rule{\imagewidth}{\imageheight}%
    }%
}
\newcommand\zoombox[2][]{
    \begin{scope}[zoombox paths]
        \pgfmathsetmacro\xpos{
            (\columncount-1)*(\imagewidth / \pgfkeysvalueof{/tikz/zoomboxarray columns} + \pgfkeysvalueof{/tikz/zoomboxarray inner gap} / \pgfkeysvalueof{/tikz/zoomboxarray columns} ) + \pgflinewidth
        }
        \pgfmathsetmacro\ypos{
            (\rowcount-1)*( \imageheight / \pgfkeysvalueof{/tikz/zoomboxarray rows} + \pgfkeysvalueof{/tikz/zoomboxarray inner gap} / \pgfkeysvalueof{/tikz/zoomboxarray rows} ) + 0.5*\pgflinewidth
        }
        \edef\dospy{\noexpand\spy [
            #1,
            zoombox paths/.append style={
                black and white pattern=\patternnumber
            },
            every spy on node/.append style={#1},
            x=\imagewidth,
            y=\imageheight
        ] on (#2) in node [anchor=north west] at ($(zoomboxes container.north west)+(\xpos pt,-\ypos pt)$);}
        \dospy
        \pgfmathtruncatemacro\pgfmathresult{ifthenelse(\columncount==\pgfkeysvalueof{/tikz/zoomboxarray columns},\rowcount+1,\rowcount)}
        \global\let\rowcount=\pgfmathresult
        \pgfmathtruncatemacro\pgfmathresult{ifthenelse(\columncount==\pgfkeysvalueof{/tikz/zoomboxarray columns},1,\columncount+1)}
        \global\let\columncount=\pgfmathresult
        \ifblackandwhitecycle
            \pgfmathtruncatemacro{\newpatternnumber}{\patternnumber+1}
            \global\edef\patternnumber{\newpatternnumber}
        \fi
    \end{scope}
}

\newtheorem{theorem}{Theorem}
\newtheorem{pr}{Proposition}
\newtheorem{lm}{Lemma}
\newtheorem{al}{Algorithm}
\newtheorem{df}{Definition}
\newtheorem{Rm}{Remark}

\begin{document}

\newcommand{\proof}{{\bf Proof:~}}
\newcommand{\uio}{$\Sigma^{(0)}$}
\newcommand{\uioo}{$\Sigma^{(1)}$}
\newcommand{\uiok}{$\Sigma^{(k)}$}
\newcommand{\eorc}{$EORC$}
\newcommand{\uiokm}{$\Sigma^{(\overline{k})}$}
\newcommand{\uioK}{$\Sigma^{(k+1)}$}

\title{Nonlinear Unknown Input Observability: The Analytic Solution in the case of a Single Unknown Input} 
\author{Agostino Martinelli
\thanks{A. Martinelli is with INRIA Rhone Alpes,
Montbonnot, France e-mail: {\tt agostino.martinelli@ieee.org}} }

\maketitle


\thanks{If you have read this book and wish to be included in a mailing list (or other means of communication) that I maintain on the subject, then send e-mail To: agostino.martinelli@inria.fr}

\tableofcontents

\begin{abstract}
This work presents the analytic solution of a fundamental open problem in the framework of nonlinear observability, which is the unknown input observability problem (UIO problem). The solution here provided holds in the case of a single unknown input.
The first part of the work presents this analytic solution, namely the analytic criterion that allows us to obtain the observability of a nonlinear system in presence of a single unknown input and multiple known inputs.
As for the observability rank condition, the proposed analytic criterion is based on the computation of the observable codistribution. Similarly to the case of only known inputs, the observable codistribution is obtained by recursively computing the Lie derivatives of the outputs along the vector fields that characterize the dynamics. However, in correspondence of the unknown input, the corresponding vector field must be suitably rescaled. Additionally, the Lie derivatives of the outputs must also be computed along a new set of vector fields that are obtained by recursively performing suitable Lie bracketing of the vector fields that define the dynamics. In practice, the entire observable codistribution is obtained by a very simple recursive algorithm. Finally, it is shown that the recursive algorithm converges in a finite number of steps and the criterion to establish that the convergence has been reached is provided. The proposed analytic extension of the observability rank condition is illustrated by checking the weak local observability of several nonlinear systems driven by multiple known inputs and a single unknown input. For these systems, extensive simulations are also provided. It is shown that, a simple estimator based on an Extended Kalman Filter, provides results that agree with what we could expect from the observability analysis. 
Finally, it is shown that the analytic criterion introduced in this work is not only a key analytical tool to automatically obtain the state observability (and so the analytic solution of a fundamental open problem) but it can be also a fundamental tool to design an appropriate estimator.

\end{abstract}

\chapter{Introduction}

State observability is a necessary condition that a state must satisfy to be estimated. Checking the system observability is a fundamental step in state estimation. 
It has become praxis in many application fields to provide an observability analysis prior to solving an estimation problem (e.g., in robotics \cite{Bry08,Hua08,Hua10,Hua13,Kwa06,Pere09,Trawny11,Vidal07,Wang08}, in visual-inertial sensor fusion \cite{Hesch12,Hesch14a,Hesch14b,Jon11,Kelly11,Kottas12,Li12,TRO12,FnT14,Mirza08,Pana13}, in sensor calibration \cite{Censi13,Guo13,Li14,ICRA06,Mirza12}). Investigating the observability properties is very simple in the linear case. Unfortunately, real systems are very rarely characterized by linearity.

The control theory community has developed the analytic tool necessary to check the state observability for nonlinear systems provided that all the system inputs are known. This is the observability rank condition introduced by Herman and Krener in 1977 \cite{Her77}. In accordance with this criterion, it is possible to obtain all the observability properties of a nonlinear system by performing automatic computation.
On the other hand, in many real scenarios, one or more disturbances can dramatically impact the system dynamics. A disturbance can be considered as an unknown input. 
Its presence can dramatically affect the observability properties of the state. This is for instance the case of a drone that operates in presence of wind. The wind is in general unknown, time-variant and acts on the system dynamics as an unknown input.

The problem of finding an analytic tool able to determine the observability properties in presence of unknown inputs is known as the {\it Unknown Input Observability} (UIO) problem. This problem was introduced and firstly investigated by the automatic control community in the seventies \cite{Basile69,Bha78,Guido71,Wang75}.  In particular, Basile and Marro provided the solution of this problem in the linear case \cite{Basile69}.
In many application fields, most of the systems are characterized by nonlinear dynamics, even in very simple cases (e.g. in planar robotics). Additionally, the presence of disturbances cannot be ignored in many cases and can significantly affect the observability properties. 
The UIO problem in the nonlinear case has recently been investigated by the robotics community and partial solutions have recently been proposed  \cite{Belo10,FnT14,ICRA14,ICRA15}. These solutions only provide sufficient conditions for the state observability. They are based on a suitable state extension.

This work provides the analytic solution of the UIO problem in the nonlinear case in the case of a single unknown input. This solution has been presented at  the SIAM symposium on control and applications 2015 \cite{SIAM-CT15}. In chapter \ref{ChapterSystemsDefinition} we define the class of systems for which we provide the solution. This class includes any nonlinear system linear in the inputs (both known and unknown).

In \cite{Her77,Isi95} the observability properties of a nonlinear system are obtained by computing the observable codistribution. The computation of this codistribution is the core of the observability rank condition introduced in \cite{Her77}. In order to deal with the case of unknown inputs, we need to derive a new algorithm able to generate the observable codistribution. In chapter \ref{ChapterObsCod}  we remind the reader the algorithm to compute the observable codistribution in the case without disturbances, together with some basic properties that characterize its convergence (section \ref{SectionORCCod}). Then, in section \ref{SectionEORCCodA},  we introduce the new algorithm that generates the entire observable codistribution in presence of a single unknown input, together with some basic properties that characterize its convergence.
The solution of the UIO problem in the case of a single unknown input is summarized in chapter \ref{ChapterEORC} and it is illustrated in chapter \ref{ChapterApplications} by checking the observability of several robotics systems driven by multiple known inputs and a single unknown input, ranging from planar robotics up to aerial navigation systems. This chapter also contains 
extensive simulations and introduces a simple estimator. The results obtained by estimating the state with this estimator fully agree with the observability analysis carried out by using the proposed analytic criterion.
In chapter \ref{ChapterProofs} we provide all the analytical derivations necessary to prove the validity of the analytic results presented in chapter \ref{ChapterObsCod}.
In chapter \ref{ChapterConclusion} we provide our conclusion by remarking that, the analytic criterion here introduced, is not only a powerful and simple criterion to automatically detect the observable states for a general nonlinear system affected by a disturbance (i.e., the analytic solution of the UIO problem in the case of a single unknown input) but it can be also a key tool for designing an appropriate estimator.


\chapter{System definition}\label{ChapterSystemsDefinition}

We will refer to a nonlinear control system with $m_u$ known inputs ($ u \triangleq [u_1,\cdots,u_{m_u}]^T$) and a single  unknown input or disturbance ($ w $). The state is the vector $ x \in M$, with $M$ an open set of $\mathbb{R}^n$. We assume that the dynamics are nonlinear with respect to the state and affine with respect to the inputs (both known and unknown). Finally, for the sake of simplicity, we will  refer to the case of a single output $y$ and we provide the extension to multiple outputs, which is trivial, separately. This will allow us to avoid the introduction of a further index. Our system is characterized by the following equations:

\begin{equation}\label{EquationSystemDefinitionDynamics1}
\left\{\begin{aligned}
  \dot{x} &=   \sum_{i=1}^{m_u}f^i ( x ) u_i +  g ( x ) w  \\
  y &= h( x ) \\
\end{aligned}\right.
\end{equation}

\noindent where  $ f^i ( x )$, $i=1,\cdots,m_u$ and $g ( x )$ are vector fields in $M$ and the function $h( x )$ is a scalar function defined on the open set $M$. Finally, we assume that the unknown input $w$ is an analytic function of time.

\noindent Throughout this document, in the case when $m_u=1$, we denote by $f(x)$ the vector field $f^1(x)$.

\chapter{Observable Codistribution}\label{ChapterObsCod}

As for the observability rank condition, the analytic criterion to obtain the observability properties of  the system in (\ref{EquationSystemDefinitionDynamics1}), is obtained by computing the observable codistribution\footnote{The reader non-familiar with the concept of {\it distribution}, as it is used in \cite{Isi95}, should not be afraid by the term {\it distribution} and the term {\it codistribution}. Very simply speaking, a distribution is a vector space defined on $M$ (our set in $\mathbb{R}^n$ where the system is defined). In particular, this vector space changes by moving on $M$. This vector space is in fact the span of a set of vector functions (vector fields) defined on $M$. A codistribution is the dual of a distribution. Very simply speaking (and this is enough to understand this manuscript) a distribution is generated by a set of column vectors. A codistribution is generated by a set of line vectors. All these vectors are vector functions (i.e., they depend on the point $x\in M$) and they have the same dimension of $x$.}. 
In this chapter we provide the algorithm to compute this codistribution. For the clarity sake, we start by reminding the reader the standard algorithm that generates the observable codisribution in absence of unknown inputs (section \ref{SectionORCCod}). Then, in section \ref{SectionEORCCodA}, we provide the algorithm for the system in (\ref{EquationSystemDefinitionDynamics1}). 

Note that this chapter directly provides the analytic results. All the analytic proofs will be provided later, in chapter \ref{ChapterProofs}.

\newpage

\section{Observable codistribution in the case without unknown inputs}\label{SectionORCCod}
We consider the system in (\ref{EquationSystemDefinitionDynamics1}) when all the inputs are known. We will denote with the symbol $\mathcal{D}$ the differential with respect to the state $x$. 
For instance, if $x=[x_1, ~x_2]^T$ and $h=x_1+x_2^2$, we have\footnote{The span of the differentials of a set of scalar functions is a codistribution. The reader non familiar with the theory of distributions can simply consider the differential as the gradient operator. The gradient of a scalar function is a line vector. For instance, if $x=[x_1, ~x_2]^T$ and $h=x_1+x_2^2$, we obtain for its gradient the line vector function $[1, ~2x_2]$. Later, in chapter \ref{ChapterApplications}, we adopt this representation. According to this, a codistribuion will be the span of a set of line vectors and a covector (i.e., an element of a codistribution) will be a line vector.}: $\mathcal{D}h=\mathcal{D}x_1+2x_2\mathcal{D}x_2$.

For a given codistribution  $\Omega$ and a given vector field $\theta=\theta(x)$ (both defined on the open set $M$), we denote by $\mathcal{L}_{\theta} \Omega$ the codistribution whose covectors are the Lie derivatives along $\theta$ of the covectors in $\Omega$. We remind the reader that the Lie derivative of a scalar function $h(x)$ along the vector field $f(x)$ is defined as follows:

\[
\mathcal{L}_f h \triangleq \frac{\partial h}{\partial x}  f
\]

\noindent which is the product of the row vector $\frac{\partial h}{\partial x} $ with the column vector $f$. Hence, it is a scalar function. Additionally, by definition of Lie derivative of covectors, we have: $\mathcal{L}_f \mathcal{D}h= \mathcal{D} \mathcal{L}_f h$.

Finally, given two vector spaces $V_1$ and $V_2$, we denote by $V_1+V_2$ their sum, i.e., the span of all the generators of both $V_1$ and $V_2$.

The observable codistribution is generated by the following recursive algorithm (see \cite{Her77} and \cite{Isi95}):

\begin{al}\label{AlgoHK}{\bf Observable codistribution when $w$ is known}

\begin{enumerate}
\item $\Omega_0=span\{\mathcal{D}h\}$;
\item $\Omega_m=\Omega_{m-1}+\sum_{i=1}^{m_u} \mathcal{L}_{f^i} \Omega_{m-1} +\mathcal{L}_{g} \Omega_{m-1} $
\end{enumerate}
\end{al}

\noindent In presence of multiple outputs, we only need to add to the codistribution $\Omega_0$, the span of the differentials of the remaining outputs.
In \cite{Isi95} it is proven that this algorithm converges. In particular, it is proven that it has converged when $\Omega_m=\Omega_{m-1}$. From this, it is easy to realize that the convergence is achieved in at most $n-1$ steps\footnote{ This is a consequence of
lemmas 1.9.1, 1.9.2 and 1.9.6 in \cite{Isi95}.}.

\section{Observable codistribution in presence of a single unknown input}\label{SectionEORCCodA}

We now provide the new algorithm that generates the observable codistribution when $w$ is an unknown input.
We will denote by $L^1_g=L^1_g(x)$ the first order Lie derivative of the function $h(x)$ along the vector field $g(x)$, i.e., 

\begin{equation}\label{EquationL1g}
L^1_g\triangleq\mathcal{L}_gh
\end{equation}

\noindent The analytic computation of the observable codistribution is based on the assumption that $L^1_g \neq 0$ on a given neighbourhood of $x_0$. In appendix \ref{AppendixCanonic}, we deal with the case when this condition is not satisfied. Specifically, we show that it is either possible to redefine the output, without altering the system observability properties, such that the Lie derivative of the new output along $g$ does not vanish or the unknown input is spurious (i.e., it does not affect the observability properties). For these reasons, we can assume that $L^1_g\neq 0$.

Before introducing the new algorithm that generates the entire observable codistribution, we introduce a new set of vector fields $^i\phi_m\in \mathbb{R}^n$ ($i=1,\cdots,m_u$ and for any integer $m$). They are obtained recursively by the following algorithm:

\begin{al}\label{AlgoPhi1}
$~$
\begin{enumerate}
\item $^i\phi_0=f^i$;
\item $^i\phi_m=\frac{[^i\phi_{m-1}, ~g]}{L^1_g}$
\end{enumerate}
\end{al}

\noindent where the parenthesis $[\cdot, \cdot]$ denote the Lie bracket of vector fields, defined as follows:

\[
[a,~b] \triangleq \frac{\partial b}{\partial x} a(x) - \frac{\partial a}{\partial x} b(x)
\]

\noindent In other words, for each $i=1,\cdots,m_u$, we have one new vector field at every step of the algorithm. Throughout this document, in the case when $m_u=1$, we denote by $\phi_m$ the vector field $^1\phi_m$.

We are now ready to provide the algorithm that generates the entire observable codistribution. It is the following:

\begin{al}\label{AlgoO1}{\bf Observable codistribution when $w$ is unknown}

\begin{enumerate}
\item $\Omega_0=span\{\mathcal{D}h\}$;
\item $\Omega_m=\Omega_{m-1}+\sum_{i=1}^{m_u} \mathcal{L}_{f^i} \Omega_{m-1} + \mathcal{L}_{\frac{g}{L^1_g}} \Omega_{m-1} +\sum_{i=1}^{m_u}span \left\{\mathcal{L}_{^i\phi_{m-1}} \mathcal{D}h\right\}$
\end{enumerate}
\end{al}

\noindent In presence of multiple outputs, we only need to add to the codistribution $\Omega_0$, the span of the differentials of the remaining outputs. Note that, in presence of multiple outputs, the function $L^1_g$ is still a scalar function since it is still defined by using a single output. The result is independent of the chosen output (provided that $L^1_g$ does not vanish\footnote{The case when $L^1_g=0$ for all the outputs is dealt in appendix \ref{AppendixCanonic}.}).

%

In section \ref{SubSectionConvergence} we investigate the convergence properties of algorithm \ref{AlgoO1}. We consider first the case of a single known input (i.e., $m_u=1$) and then the results are easily extended to the case of multiple known inputs ($m_u>1$) in section \ref{SubSectionExt}.
We prove that algorithm \ref{AlgoO1} converges and we also provide the analytic criterion to check that the convergence has been attained.
This proof and the convergence criterion cannot be the same that hold for algorithm \ref{AlgoHK}, because of the last term that appears at the recursive step\footnote{The convergence criterion of algorithm \ref{AlgoHK} is a consequence of the fact that, all the terms that appear at the recursive step of algorithm \ref{AlgoHK}, are the Lie derivative of the codistribution at the previous step, along fixed vector fields (i.e., vector fields that remain the same at each step of the algorithm). This is not the case for the last term at the recursive step of algorithm \ref{AlgoO1}.}, i.e., the term $\sum_{i=1}^{m_u}span \left\{\mathcal{L}_{^i\phi_{m-1}} \mathcal{D}h\right\}$ (the special case when, the contribution due to this last term is included in the other terms, is considered separately by lemma \ref{LemmaConvergenceSpecial}).
In general, the criterion to establish that the convergence has been attained is not simply obtained by checking the condition $\Omega_{m+1}=\Omega_m$.
Deriving the new analytic criterion is not immediate. It requires to derive the analytic expression that describes 
the behaviour of the last term in the recursive step. 
This fundamental equation is provided in chapter \ref{ChapterProofs} and it is the equation (\ref{EquationKeyEquality}). 
The analytic derivation of this equation allows us to detect the key quantity that governs the convergence of algorithm \ref{AlgoO1}, in particular regarding the contribution due to the last term in the recursive step. This key quantity is the following scalar:

\begin{equation}\label{EquationRho}
\tau=\frac{\mathcal{L}^2_gh}{(L^1_g)^2}~~
\end{equation}

\noindent  We prove (see lemma \ref{LemmaRhoInOmegam} in chapter \ref{ChapterProofs}) that, in general, it exists $m'$ such that $\mathcal{D}\tau \in \Omega_{m'}$ (and therefore $\mathcal{D}\tau \in \Omega_m$ $\forall m\ge m'$). Additionally, we prove that the convergence of the algorithm has been reached when $\Omega_{m+1}=\Omega_m$, $m\ge m'$ and $m\ge 2$ (theorem \ref{TheoremStop}). We also prove that the required number of steps is at most $n+2$.

In section \ref{SubSectionSeparation} it is also shown that the computed codistribution is the entire observable codistribution. Also in this case, the proof is given by first considering the case of a single known input (see theorem \ref{TheoremSeparation}) and then, its validity is extended to the case of multiple inputs in section \ref{SubSectionExt}. Note that this proof is based on the assumption that the unknown input ($w$) is a differentiable function of time, up to a given order (the order depends on the specific case).

Algorithm \ref{AlgoO1} differs from the standard algorithm \ref{AlgoHK} because of the following reasons:

\begin{itemize}

\item In the recursive step, the vector field that corresponds to the unknown input (i.e., $g$) must be rescaled by dividing by $L^1_g$.

\item The recursive step also contains the sum of the contributions $\sum_{i=1}^{m_u}\mathcal{L}_{^i\phi_{m-1}} \mathcal{D}h$. In other words, we need to compute the Lie derivatives of the differential of the output along the vector fields obtained through the recursive algorithm \ref{AlgoPhi1}.

\item The convergence of algorithm \ref{AlgoO1} is achieved in at most $n+2$ steps, instead of $n-1$ steps (in the special case dealt by lemma \ref{LemmaConvergenceSpecial}, this upper bound is $n-1$ for both cases).

\item When $\Omega_m=\Omega_{m-1}$ algorithm \ref{AlgoHK} has converged. For algorithm \ref{AlgoO1}, we also need to check that $\mathcal{D} \tau \in \Omega_m$ and $m\ge 2$ (with the exception of the special case dealt by lemma \ref{LemmaConvergenceSpecial}).

\end{itemize}

\chapter{The Analytic Criterion}\label{ChapterEORC}

%

In this chapter, we outline all the steps to investigate the weak local observability at a given point $x_0$ of a nonlinear system characterized by (\ref{EquationSystemDefinitionDynamics1}). Basically, these steps are the steps necessary to compute the observable codistribution (i.e., the steps of algorithm \ref{AlgoO1}) and to prove that the differential of a given state component belongs to this codistribution.

Note that, in the trivial case analyzed by lemma \ref{LemmaConvergenceSpecial}, the criterion provided below simplifies, since we do not need to compute the quantity $\tau\left(=\frac{\mathcal{L}^2_g h}{(\mathcal{L}^1_g h)^2}\right)$, and we do not need to check that its differential belongs to the codistribution computed at every step of algorithm \ref{AlgoO1}. In practice, we skip the steps $4$ and $5$ in the procedure below.

\begin{enumerate}

\item For the chosen $x_0$, compute $L^1_g(= \mathcal{L}^1_g h)$. In the case when $L^1_g=0$, choose another function in the space of functions $\mathcal{F}$ (defined as the space that contains $h$ and its Lie derivative up to any order along the vector fields $f^1,\cdots,f^{m_u}$) such that its Lie derivative along $g$ does not vanish\footnote{If the Lie derivative of any function in $\mathcal{F}$ vanishes,  it means that the unknown input can be ignored to obtain the observability properties, as it is shown in appendix \ref{AppendixCanonic}.}.

\item  Compute the codistribution $\Omega_0$ and $\Omega_1$ (at $x_0$) by using algorithm \ref{AlgoO1}.

\item  Compute the vector fields $^i\phi_m$ ($i=1,\cdots,m_u$) by using algorithm \ref{AlgoPhi1}, starting from $m=0$, to check if the considered system is in the special case dealt by lemma \ref{LemmaConvergenceSpecial}\footnote{Note that this check requires at most $n-1$ (simple) operations.}. In this trivial case, set $m'=0$,  use the recursive step of algorithm \ref{AlgoO1} to build the codistribution $\Omega_m$ for $m\ge 2$, and skip to step $6$.

\item Compute $\tau\left(=\frac{\mathcal{L}^2_g h}{(L^1_g)^2}\right)$ and $\mathcal{D} \tau$.

\item Use the recursive step of algorithm \ref{AlgoO1} to build the codistribution $\Omega_m$ for $m\ge 2$, 
and, for each $m$, check if $\mathcal{D}\tau \in \Omega_m$. Denote by $m'$ the smallest $m$ such that $\mathcal{D}\tau \in \Omega_m$.

\item For each $m\ge m'$, check if $\Omega_{m+1}=\Omega_m$ and denote by $\Omega^*=\Omega_{m^*}$, where $m^*$ is the smallest integer such that $m^*\ge m'$ and $\Omega_{m^*+1}=\Omega_{m^*}$ (note that $m^*\le n+2$).

\item If the differential of a given state component ($x_j$, $j=1,\cdots,n$) belongs to $\Omega^*$ (namely if $\mathcal{D}x_j\in\Omega^*$) on a given neighbourhood of $x_0$, then $x_j$ is weakly locally observable at $x_0$. If this holds for all the state components, the state $x$ is weakly locally observable at $x_0$. Finally, if the dimension of $\Omega^*$ is smaller than $n$ on a given neighbourhood of $x_0$, then the state is not weakly locally observable at $x_0$.

\end{enumerate}

\chapter{Applications}\label{ChapterApplications}

We apply the criterion described in chapter \ref{ChapterEORC} in order to investigate the observability properties of several nonlinear systems characterized by the equations given in (\ref{EquationSystemDefinitionDynamics1}). 
These examples will be discussed in sections \ref{SectionApplicationUnicycleUI}, \ref{SectionApplicationUnicycleSD} and \ref{SectionApplicationDefinition}. The dynamics of  the first two examples are the ones of the unicycle. In accordance with the unicycle dynamics, the motion is powered by two independent controls, which are the linear and the angular speed, respectively. 
In section \ref{SectionApplicationUnicycleUI}
we consider the case when one of these two inputs is unknown and acts as an unknown input.
In section \ref{SectionApplicationUnicycleSD}
we consider the case when both these two inputs are known. However, the dynamics are also affected by an external unknown input.
Finally, the last example regards a vehicle that moves in $3D$ (section \ref{SectionApplicationDefinition}).

These  examples are deliberately uncomplicated in order to allow us to compare the analytic results provided by the proposed criterion with what we can expect by following intuitive reasoning.

Note that the analytic criterion is very powerful and can be used to automatically obtain the observability properties of any system that satisfies (\ref{EquationSystemDefinitionDynamics1}), i.e., independently of the state dimension (intuitive reasoning becomes often prohibitive for systems characterized by high-dimensional states).

This chapter also contains 
extensive simulations and introduces a simple estimator. The results obtained by estimating the state with this estimator fully agree with the observability analysis carried out by using the proposed analytic criterion. These results  also show that the analytic criterion is not only a criterion to check the state observability but also a fundamental tool to design an appropriate estimator for a given nonlinear system in presence of a disturbance (see chapter \ref{ChapterConclusion} for a discussion about this fundamental aspect).

\newpage

\section{Unicycle with one input unknown}\label{SectionApplicationUnicycleUI}

\subsection{The system}\label{SubSectionApplicationUnicycleUI}

We consider a vehicle that moves on a $2D$-environment. The configuration of the vehicle in a global reference frame, can be characterized through the vector $[x_v, ~y_v, ~\theta]^T$ where $x_v$ and $y_v$ are the Cartesian vehicle coordinates, and $\theta$ is the vehicle orientation. We assume that the dynamics of this vector satisfy the unicycle differential equations:

\begin{equation}
\left[\begin{aligned}
  \dot{x}_v ~&= v \cos\theta \\
  \dot{y}_v ~&= v \sin\theta \\
  \dot{\theta} ~~&= \omega \\
\end{aligned}\right.
\end{equation}

\noindent where $v$ and $\omega$ are the linear and the rotational vehicle speed, respectively, and they are the system inputs. We consider the following three cases of output (see also figure \ref{FigUI} for an illustration):

\begin{figure}[htbp]
\begin{center}
\includegraphics[width=.8\columnwidth]{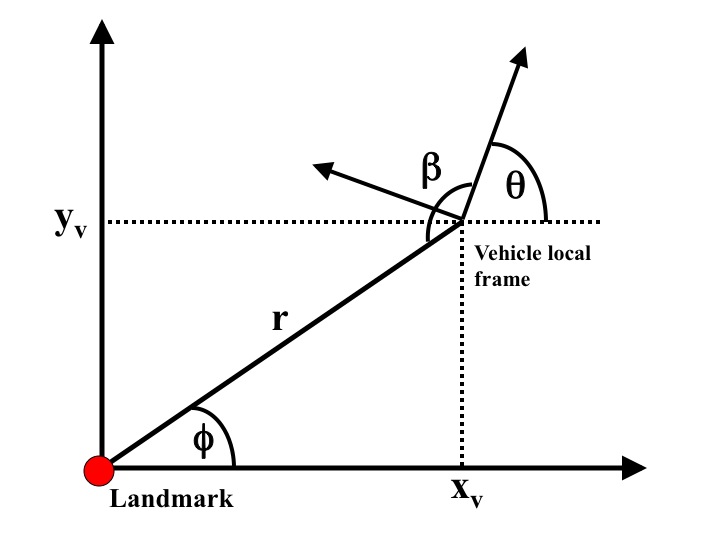}
\caption{The vehicle state in cartesian and polar coordinates ($[x_v,y_v,\theta]^T$ and $[r,\phi,\theta]^T$, respectively) together with the three considered outputs ($r,\beta$ and $\phi$).} \label{FigUI}
\end{center}
\end{figure}

\begin{enumerate}

\item the distance ($r$) from the origin (e.g., a landmark is at the origin and its distance is measured by a range sensor);

\item the bearing angle ($\beta$) of the origin in the local frame (e.g., a landmark is at the origin and its bearing angle is measured by an on-board camera);

\item the bearing angle ($\phi$) of the vehicle in the global frame (e.g., a camera is placed at the origin).

\end{enumerate}

\noindent We can analytically express the output in terms of the state. We remark that the expressions become very simple if we adopt polar coordinates: $r \triangleq \sqrt{x_v^2+y_v^2}$, $\phi=atan\frac{y_v}{x_v}$. We have, for the three cases, $y=r$, $y=\pi-(\theta-\phi)$ and $y=\phi$, respectively. For each of these three cases, we consider the following two cases: $v$ is known, $\omega$ is unknown;  $v$ is unknown, $\omega$ is known. Hence, we have six cases. The dynamics of the unicycle in polar coordinates become:

\begin{equation}\label{EquationSimpleExampeDynamics}
\left[\begin{aligned}
  \dot{r} &= v \cos(\theta-\phi) \\
  \dot{\phi} &= \frac{v}{r} \sin(\theta-\phi) \\
  \dot{\theta} &= \omega \\
\end{aligned}\right.
\end{equation}

\subsection{Intuitive procedure to obtain the observability properties}\label{SubSectionIntuitive}

By using the observability rank condition in \cite{Her77}, we easily obtain that, when both the inputs are known, the dimension of the observable codistribution is $2$ for the first two observations ($y=r$ and $y=\theta-\phi$) and $3$ for the last one ($y=\phi$). In particular, for the first two observations all the initial states rotated around the vertical axis are indistinguishable. When one of the inputs misses, this unobservable degree of freedom obviously remains. On the other hand, when the linear speed is unknown (i.e., it acts as an unknown input ($w=v$)) and the observation is an angle (second and third observation, i.e., $y=\theta-\phi$ and $y=\phi$, respectively), we lose a further degree of freedom, which corresponds to the absolute scale. In table \ref{Table1} we provide the dimension of the observable codistribuion obtained by following this intuitive reasoning for the six considered cases.

\begin{table}
\begin{center}
\begin{tabular}{|l|c|}
  \hline
  Case & Dimension of $\Omega$\\
  \hline
  \hline
  $1^{st}$: $y=r$, $u=\omega$, $w=v$ & 2 \\
  \hline
  $2^{nd}$: $y=r$, $w=\omega$, $u=v$ & 2 \\
  \hline
  $3^d$: $y=\theta -\phi$, $u=\omega$, $w=v$ & 1 \\
  \hline
  $4^{th}$: $y=\theta -\phi$, $w=\omega$, $u=v$ & 2 \\
  \hline
  $5^{th}$: $y=\phi$, $u=\omega$, $w=v$ & 2 \\
  \hline
  $6^{th}$: $y=\phi$, $w=\omega$, $u=v$ & 3 \\
  \hline
\end{tabular}
\end{center}
\caption{Dimension of the observable codistribution ($\Omega$) obtained by following intuitive reasoning}
\label{Table1}
\end{table}

\subsection{Analytic results}\label{SubSectionAnalytic}

We now derive the observability properties by applying the analytic criterion described in chapter \ref{ChapterEORC}. For all the cases we have $m_u=1$. Hence, we adopt the following notation: $f\triangleq f^1$ (for the vector field in (\ref{EquationSystemDefinitionDynamics1})) and $\phi_m \triangleq ~^1\phi_m$ (for the vectors defined by algorithm \ref{AlgoPhi1}).
We consider the six cases defined in section \ref{SubSectionApplicationUnicycleUI}, separately.

\vskip.5cm
\subsubsection{First Case: $y=r$, $u=\omega$, $w=v$} 

We have:

\[
f=\left[\begin{array}{c}
 0 \\
 0  \\
  1 \\
\end{array}
\right]~~~g=\left[\begin{array}{c}
 \cos(\theta-\phi) \\
 \frac{\sin(\theta-\phi)}{r}  \\
  0 \\
\end{array}
\right]
\]

\noindent We apply the analytic criterion in chapter \ref{ChapterEORC}. We obtain:

\subsubsection{Step 1}

We have $L^1_g=\cos(\theta-\phi)$, which does not vanish, in general.

\subsubsection{Step 2}

We have: $\Omega_0=span\{[1,0,0]\}$. Additionally, $\Omega_1=\Omega_0$.

\subsubsection{Step 3}

We have $\mathcal{L}_{\phi_0}L^1_g=\mathcal{L}_fL^1_g=-\sin(\theta-\phi)$, which does not vanish, in general. This suffices to conclude that the considered system is not in the special case considered by lemma \ref{LemmaConvergenceSpecial} and we need to continue with step 4.

\subsubsection{Step 4}

We have
$\tau\triangleq \frac{L^2_g}{(L^1_g)^2}=\frac{\tan^2(\theta-\phi)}{r}$
and
\[
\mathcal{D}\tau=\frac{\tan(\theta-\phi)}{r}\left[-\frac{\tan(\theta-\phi)}{r}, -\frac{2}{\cos^2(\theta-\phi)},\frac{2}{\cos^2(\theta-\phi)}\right]
\]

\subsubsection{Step 5}

\noindent We need to compute $\Omega_2$ and, in order to do this, we need to compute $\phi_1$. We obtain:

\[
\phi_1=\left[\begin{array}{c}
  -\tan(\theta-\phi) \\
  \frac{1}{r}  \\
  0 \\
\end{array}
\right]
\]

\noindent and 

\[
\Omega_2=span\left\{[1,0,0], ~\left[0,\frac{1}{\cos^2(\theta-\phi)}, -\frac{1}{\cos^2(\theta-\phi)}\right]\right\}
\]

\noindent It is immediate to check that  $\mathcal{D}\tau\in \Omega_2$, meaning that $m'=2$.

\subsubsection{Step 6}

By a direct computation, it is possible to check that $\Omega_3=\Omega_2$ meaning that $m^*=2$ and $\Omega^*=\Omega_2$

\subsubsection{Step 7}

The dimension of the observable codistribution is $2$. We conclude that the state is not weakly locally observable. This result agrees with the one in table \ref{Table1} (second line).

\vskip.2cm
\subsubsection{Second Case: $y=r$, $u=v$, $w=\omega$} 

We have:

\[
f=\left[\begin{array}{c}
 \cos(\theta-\phi) \\
 \frac{\sin(\theta-\phi)}{r}  \\
  0 \\
\end{array}
\right]~~~g=\left[\begin{array}{c}
 0 \\
 0  \\
  1 \\
\end{array}
\right]
\]

\noindent We apply the analytic criterion in chapter \ref{ChapterEORC}. We obtain:

\subsubsection{Step 1}

We easily obtain $\mathcal{L}_g h=0$.
We consider the function $\mathcal{L}_fh\in \mathcal{F}$. We have: $\mathcal{L}_fh=\cos(\theta-\phi)$ and $\mathcal{L}_g\mathcal{L}_fh=-\sin(\theta-\phi)$, which does not vanish, in general. Hence, we can proceed with the steps in chapter \ref{ChapterEORC} by setting:

\[
h=\cos(\theta-\phi)
\]

\noindent We obtain: $L^1_g=-\sin(\theta-\phi)$

\subsubsection{Step 2}

We have:

\[
\Omega_0=span\{[1, 0, 0],~ [0, \sin(\theta-\phi), ~-\sin(\theta-\phi)]\}
\]

\noindent as long as the function $r$ is also a system output. Additionally, $\Omega_1=\Omega_0$.

\subsubsection{Step 3}

We have $\mathcal{L}_{\phi_0}L^1_g=\mathcal{L}_fL^1_g=\frac{\sin(\theta-\phi)\cos(\theta-\phi)}{r}$, which does not vanish, in general. This suffices to conclude that the considered system is not in the special case considered by lemma \ref{LemmaConvergenceSpecial} and we need to continue with step 4.

\subsubsection{Step 4}

We have:

\[
\tau=-\frac{\cos(\theta-\phi)}{\sin^2(\theta-\phi)}
\]

\subsubsection{Step 5}

\noindent By a direct computation we obtain 
$\Omega_2=\Omega_1$. Additionally, it is immediate to check that  $\mathcal{D}\tau\in \Omega_2$, meaning that $m'=2$. 

\subsubsection{Step 6}

By a direct computation, it is possible to check that $\Omega_3=\Omega_2$ meaning that $m^*=2$ and $\Omega^*=\Omega_2$.

\subsubsection{Step 7}

The dimension of the observable codistribution is $2$. We conclude that the state is not weakly locally observable. This result agrees with the one in table \ref{Table1} (third line).

\vskip.2cm
\subsubsection{Third Case: $y=\theta-\phi$, $u=\omega$, $w=v$} 
We have:

\[
f=\left[\begin{array}{c}
 0 \\
 0  \\
  1 \\
\end{array}
\right]~~~g=\left[\begin{array}{c}
 \cos(\theta-\phi) \\
 \frac{\sin(\theta-\phi)}{r}  \\
  0 \\
\end{array}
\right]
\]

\noindent We apply the analytic criterion in chapter \ref{ChapterEORC}. We obtain:

\subsubsection{Step 1}

We have $L^1_g=-\frac{\sin(\theta-\phi)}{r}$, which does not vanish, in general.

\subsubsection{Step 2}

We have $\Omega_0=span\{[0,-1,1]\}$ and $\Omega_1=\Omega_0$.

\subsubsection{Step 3}

We have $\mathcal{L}_{\phi_0}L^1_g=\mathcal{L}_fL^1_g=-\frac{\cos(\theta-\phi)}{r}$, which does not vanish, in general. This suffices to conclude that the considered system is not in the special case considered by lemma \ref{LemmaConvergenceSpecial} and we need to continue with step 4.

\subsubsection{Step 4}

We have $\tau=2\cot(\theta-\phi)$ and 

\[
\mathcal{D}\tau=\frac{2}{\sin^2(\theta-\phi)}\left[0, 1,-1\right]
\]

\subsubsection{Step 5}

\noindent By a direct computation we obtain 
$\Omega_2=\Omega_1$. Additionally, it is immediate to check that  $\mathcal{D}\tau\in \Omega_2$, meaning that $m'=2$. 

\subsubsection{Step 6}

By a direct computation, it is possible to check that $\Omega_3=\Omega_2$ meaning that $m^*=2$ and $\Omega^*=\Omega_2$.

\subsubsection{Step 7}

The dimension of the observable codistribution is $1$. We conclude that the state is not weakly locally observable. This result agrees with the one in table \ref{Table1} (fourth line). Note that, the new unobservable direction with respect to the case when both inputs are known, is precisely the absolute scale, since the vector $\mathcal{D} r = [1,0,0]\notin \Omega^*$.

\vskip.2cm
\subsubsection{Fourth Case: $y=\theta-\phi$, $u=v$, $w=\omega$} We have:

\[
f=\left[\begin{array}{c}
 \cos(\theta-\phi) \\
 \frac{\sin(\theta-\phi)}{r}  \\
  0 \\
\end{array}
\right]~~~g=\left[\begin{array}{c}
 0 \\
 0  \\
  1 \\
\end{array}
\right]
\]

\noindent We apply the analytic criterion in chapter \ref{ChapterEORC}. We obtain:

\subsubsection{Step 1}

We have $L^1_g=1\neq 0$.

\subsubsection{Step 2}

By a direct computation we obtain:
$\Omega_0=span\{[0,-1,1]\}$ and

\[
\Omega_1=span\left\{[0,-1,1],  \left[-\frac{\sin(\theta-\phi)}{r^2},-\frac{\cos(\theta-\phi)}{r}, \frac{\cos(\theta-\phi)}{r}\right]\right\}
\]

\subsubsection{Step 3}

Since $L^1_g=1$, it is immediate to realize that $\mathcal{L}_{\phi_j}L^1_g=0$, for any integer $j\ge 0$. Hence, the considered system is in the special case considered by lemma \ref{LemmaConvergenceSpecial} and we skip to step 6 by setting $m'=0$.

\subsubsection{Step 6}

By a direct computation, it is possible to check that $\Omega_2=\Omega_1$ meaning that $m^*=1$ and $\Omega^*=\Omega_1$.

\subsubsection{Step 7}

The dimension of the observable codistribution is $2$. We conclude that the state is not weakly locally observable. This result agrees with the one in table \ref{Table1} (fifth line).

\vskip.2cm
\subsubsection{Fifth Case: $y=\phi$, $u=\omega$, $w=v$} We have:

\[
f=\left[\begin{array}{c}
 0 \\
 0  \\
  1 \\
\end{array}
\right]~~~g=\left[\begin{array}{c}
 \cos(\theta-\phi) \\
 \frac{\sin(\theta-\phi)}{r}  \\
  0 \\
\end{array}
\right]
\]

\noindent We apply the analytic criterion in chapter \ref{ChapterEORC}. We obtain:

\subsubsection{Step 1}

We have $L^1_g=\frac{\sin(\theta-\phi)}{r}$, which does not vanish, in general.

\subsubsection{Step 2}

We easily obtain $\Omega_0=span\{[0,1,0]\}$ and $\Omega_1=\Omega_0$.

\subsubsection{Step 3}

We have $\mathcal{L}_{\phi_0}L^1_g=\mathcal{L}_fL^1_g=\frac{\cos(\theta-\phi)}{r}$, which does not vanish, in general. This suffices to conclude that the considered system is not in the special case considered by lemma \ref{LemmaConvergenceSpecial} and we need to continue with step 4.

\subsubsection{Step 4}

We have $\tau=-2\cot(\theta-\phi)$ and

\[
\mathcal{D}\tau=\frac{2}{\sin^2(\theta-\phi)}\left[0, -1,1\right]
\]

\subsubsection{Step 5}

To compute $\Omega_2$ we need to compute $\phi_1$. We obtain:

\[
\phi_1=\left[\begin{array}{c}
  -r \\
  \cot(\theta-\phi)  \\
  0 \\
\end{array}
\right]
\]

\noindent and

\[
\Omega_2=span\left\{[0,1,0], ~\frac{1}{\sin^2(\theta-\phi)}\left[0, 1,-1\right]\right\}
\]

\noindent It is immediate to check that  $\mathcal{D}\tau\in \Omega_2$, meaning that $m'=2$. 

\subsubsection{Step 6}

By a direct computation we obtain $\Omega_3=\Omega_2$ meaning that $m^*=2$ and $\Omega^*=\Omega_2$, whose dimension is $2$.

\subsubsection{Step 7}

The dimension of the observable codistribution is $2$. We conclude that the state is not weakly locally observable. This result agrees with the one in table \ref{Table1} (sixth line). Note that, the new unobservable direction with respect to the case when both inputs are known, is precisely the absolute scale, since the vector $\mathcal{D} r = [1,0,0]\notin \Omega^*$.

\vskip.2cm
\subsubsection{Sixth Case: $y=\phi$, $u=v$, $w=\omega$} We have:

\[
f=\left[\begin{array}{c}
 \cos(\theta-\phi) \\
 \frac{\sin(\theta-\phi)}{r}  \\
  0 \\
\end{array}
\right]~~~g=\left[\begin{array}{c}
 0 \\
 0  \\
  1 \\
\end{array}
\right]
\]

\noindent We apply the analytic criterion in chapter \ref{ChapterEORC}. We obtain:

\subsubsection{Step 1}

We easily obtain $\mathcal{L}_g h=0$.
We consider the function $\mathcal{L}_fh\in \mathcal{F}$. We have: $\mathcal{L}_fh=\frac{\sin(\theta-\phi)}{r}$ and $\mathcal{L}_g\mathcal{L}_fh=\frac{\cos(\theta-\phi)}{r}$, which does not vanish, in general. Hence, we can proceed with the steps in chapter \ref{ChapterEORC} by setting:

\[
h=\frac{\sin(\theta-\phi)}{r}~~~~L^1_g=\frac{\cos(\theta-\phi)}{r}
\]

\subsubsection{Step 2}

We have:

\[
\Omega_0=span\left\{[0, 1, 0],~ \left[-\frac{\sin(\theta-\phi)}{r^2}, -\frac{\cos(\theta-\phi)}{r}, \frac{\cos(\theta-\phi)}{r}\right]\right\}
\]

\noindent as long as the function $\phi$ is also a system output. We compute $\Omega_1$. By a direct computation, we obtain that its dimension is $3$.
Hence, we do not need to proceed with the remaining steps since we can directly conclude that the entire state is weakly locally observable. This result agrees with the one in table \ref{Table1} (seventh line).


\section{Unicycle in presence of an external unknown input}\label{SectionApplicationUnicycleSD}

\subsection{The system}\label{SubSectionApplicationSingle}

\begin{figure}[htbp]
\begin{center}
\includegraphics[width=.8\columnwidth]{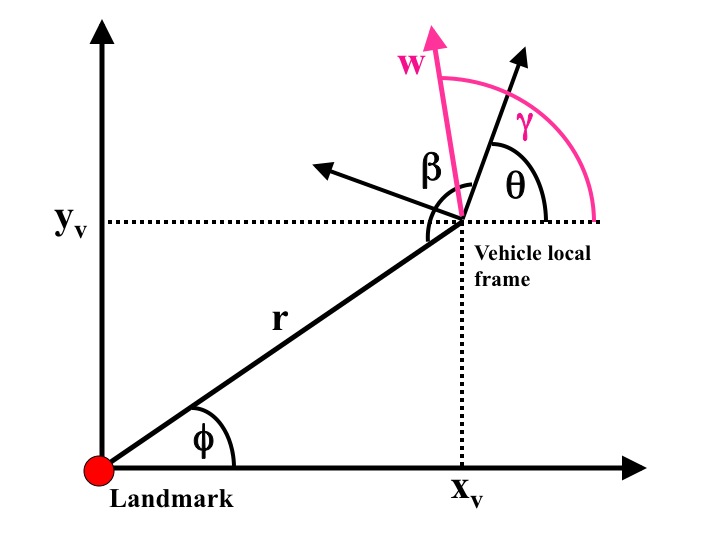}
\caption{The vehicle state together with the three considered outputs.} \label{FigSD}
\end{center}
\end{figure}

We consider the same vehicle considered in section \ref{SectionApplicationUnicycleUI} and we adopt the same state $[x_v, ~y_v, ~\theta]^T$ to characterize its position and orientation.
We assume that the vehicle motion is also affected by an unknown input that produces an additional (and unknown) robot speed (denoted by $w$) along a fixed direction (denoted by $\gamma$). Hence, the dynamics are characterized by the following differential equations:

\begin{equation}\label{EquationSystemDynamics}
\left[\begin{aligned}
  \dot{x}_v ~&= v \cos\theta +w \cos\gamma\\
  \dot{y}_v ~&= v \sin\theta +w \sin\gamma\\
  \dot{\theta} ~~&= \omega \\
\end{aligned}\right.
\end{equation}

\noindent where $v$ and $\omega$ are the linear and the rotational speed, respectively, in absence of the unknown input. We assume that these two speeds are known (we refer to them as to the known inputs), $w$ is unknown (we refer to it as to the unknown input or disturbance) and $\gamma$ is constant in time.  See also figure \ref{FigSD} for an illustration.

We consider the same three cases of output considered in section \ref{SectionApplicationUnicycleUI}. Additionally, we deal with both the cases when $\gamma$ is known and unknown (in section \ref{SubSectionObservabilityKnown} and \ref{SubSectionObservabilityUnknown}, respectively).

\subsection{Observability properties when $\gamma$ is known}\label{SubSectionObservabilityKnown}

\noindent The state  is $[x_v, ~y_v, ~\theta]^T$ and its dynamics are provided by the three equations in (\ref{EquationSystemDynamics}), where $\gamma$ is a known parameter. 
These equations are a special case of (\ref{EquationSystemDefinitionDynamics1}).
From (\ref{EquationSystemDynamics}) and (\ref{EquationSystemDefinitionDynamics1}) we easily obtain: $m_u=2$, $u_1=v$, $u_2=\omega$,

\[
f^1=\left[\begin{array}{c}
 \cos\theta \\
 \sin\theta \\
 0  \\
\end{array}
\right],~~~f^2=\left[\begin{array}{c}
 0 \\
 0  \\
  1 \\
\end{array}
\right],~~~g=\left[\begin{array}{c}
 \cos\gamma \\
 \sin\gamma \\
  0 \\
\end{array}
\right]
\]

\noindent We consider the three outputs separately. For the simplicity sake, we actually consider the following three outputs: $y=r^2=x_v^2+y_v^2$ instead of $y=r$, $y=\tan\beta=\frac{y_v-x_v\tan\theta}{x_v+y_v\tan\theta}$ instead of $y=\beta$ and $y=\tan\phi=\frac{y_v}{x_v}$ instead of $y=\phi$. Obviously, the result of the observability analysis does not change.

\vskip.2cm
\subsubsection{First Case: $y=r^2$} 
We apply the analytic criterion in chapter \ref{ChapterEORC}. We obtain:

\subsubsection{Step 1}

We have $L^1_g=2(x_v\cos\gamma + y_v\sin\gamma)$, which does not vanish, in general.

\subsubsection{Step 2}

We have $\Omega_0=span\{[x_v,y_v,0]\}$ and $\Omega_1=span\{[x_v,y_v,0], [\cos\theta, \sin\theta, y_v\cos\theta-x_v\sin\theta]~\}$.

\subsubsection{Step 3}

We have $\mathcal{L}_{^1\phi_0}L^1_g=\mathcal{L}_{f^1}L^1_g=2\cos(\gamma-\theta)$, which does not vanish, in general. This suffices to conclude that the considered system is not in the special case considered by lemma \ref{LemmaConvergenceSpecial} and we need to continue with step 4.

\subsubsection{Step 4}

We have $\tau= \frac{1}{2(x_v\cos\gamma + y_v\sin\gamma)^2}$ and 

\[
\mathcal{D}\tau=-\frac{1}{(x_v\cos\gamma + y_v\sin\gamma)^3}
\left[\cos\gamma, \sin\gamma,0\right]
\]

\subsubsection{Step 5}

We need to compute $\Omega_2$ and, in order to do this, we need to compute $^1\phi_1$ and $^2\phi_1$ through algorithm \ref{AlgoPhi1}. We obtain: 
$^1\phi_1=^2\phi_1=\left[\begin{array}{c}
  0 \\
  0  \\
  0 \\
\end{array}
\right]$. On the other hand, we obtain that $\mathcal{L}_{\frac{g}{L^1_g}}\mathcal{D} \mathcal{L}_{f^1}h \notin \Omega_1$. Hence, by using algorithm \ref{AlgoO1} we obtain that $\Omega_2$ has dimension equal to $3$. As a result, we do not need to proceed with the remaining steps, since we can directly conclude that the entire state is weakly locally observable.

\vskip.2cm
\subsubsection{Second Case: $y=\tan\beta$} 
We apply the analytic criterion in chapter \ref{ChapterEORC}. We obtain:

\subsubsection{Step 1}

We have $L^1_g=-\frac{y_v \cos\gamma - x_v \sin\gamma}{x_v^2 \cos^2\theta + 2 \sin\theta \cos\theta x_v y_v  - y_v^2 \cos^2\theta + y_v^2}$, which does not vanish, in general.

\subsubsection{Step 2}

By an explicit computation (by using algorithm \ref{AlgoO1}) we obtain that the dimension of $\Omega_0$ is $1$ and the dimension of $\Omega_1$ is $2$. 

\subsubsection{Step 3}

We have $\mathcal{L}_{^2\phi_0}L^1_g=\mathcal{L}_{f^2}L^1_g\neq 0$, in general. This suffices to conclude that the considered system is not in the special case considered by lemma \ref{LemmaConvergenceSpecial} and we need to continue with step 4.

\subsubsection{Step 4}

We have $\tau=\frac{x_v\cos\gamma + y_v\sin\gamma + x_v\cos(\gamma - 2\theta) - y_v\sin(\gamma - 2\theta)}{y_v\cos\gamma - x_v\sin\gamma}$

\subsubsection{Step 5}

We need to compute $\Omega_2$. Also in this case, we obtain that $\mathcal{L}_{\frac{g}{L^1_g}}\mathcal{D} \mathcal{L}_{f^1}h \notin \Omega_1$. Hence, by using algorithm \ref{AlgoO1} we obtain that $\Omega_2$ has dimension equal to $3$. As a result, we do not need to proceed with the remaining steps, since we can directly conclude that the entire state is weakly locally observable.

\vskip.2cm
\subsubsection{Third Case: $y=\tan\phi$} 
We apply the analytic criterion in chapter \ref{ChapterEORC}. We obtain:

\subsubsection{Step 1}

We have $L^1_g=-\frac{y_v\cos\gamma - x_v\sin\gamma}{x_v^2}$, which does not vanish, in general.

\subsubsection{Step 2}

We have $\Omega_0=span\{[-y_v,x_v,0]\}$. In addition, by an explicit computation (by using algorithm \ref{AlgoO1}) we obtain that the dimension of $\Omega_1$ is $2$. 

\subsubsection{Step 3}

We have $\mathcal{L}_{^2\phi_0}L^1_g=\mathcal{L}_{f^2}L^1_g\neq 0$, in general. This suffices to conclude that the considered system is not in the special case considered by lemma \ref{LemmaConvergenceSpecial} and we need to continue with step 4.

\subsubsection{Step 4}

We have $\tau= -\frac{2x_v^4\cos\gamma}{x_v^4\sin\gamma - x_v^3y_v\cos\gamma}$

\subsubsection{Step 5}

We need to compute $\Omega_2$. Also in this case, we obtain that $\mathcal{L}_{\frac{g}{L^1_g}}\mathcal{D} \mathcal{L}_{f^1}h \notin \Omega_1$. Hence, by using algorithm \ref{AlgoO1} we obtain that $\Omega_2$ has dimension equal to $3$. As a result, we do not need to proceed with the remaining steps, since we can directly conclude that the entire state is weakly locally observable.

\subsection{Observability properties when $\gamma$ is unknown}\label{SubSectionObservabilityUnknown}

\noindent The state  is $[x_v, ~y_v, ~\theta, ~\gamma]^T$ and its dynamics are provided by the following four equations:

\begin{equation}\label{EquationSystemDynamicsGamma}
\left[\begin{aligned}
  \dot{x}_v ~&= v \cos\theta +w \cos\gamma\\
  \dot{y}_v ~&= v \sin\theta +w \sin\gamma\\
  \dot{\theta} ~~&= \omega \\
  \dot{\gamma} ~~&= 0 \\
\end{aligned}\right.
\end{equation}

From (\ref{EquationSystemDynamicsGamma}) and (\ref{EquationSystemDefinitionDynamics1}) we easily obtain: $m_u=2$, $u_1=v$, $u_2=\omega$,

\[
f^1=\left[\begin{array}{c}
 \cos\theta \\
 \sin\theta \\
 0  \\
 0  \\
\end{array}
\right],~~~f^2=\left[\begin{array}{c}
 0 \\
 0  \\
  1 \\
 0  \\
\end{array}
\right],~~~g=\left[\begin{array}{c}
 \cos\gamma \\
 \sin\gamma \\
  0 \\
 0  \\
\end{array}
\right]
\]

\noindent We consider the three outputs separately. As in section \ref{SubSectionObservabilityKnown}, we actually consider the following three outputs: $y=r^2=x_v^2+y_v^2$ instead of $y=r$, $y=\tan\beta=\frac{y_v-x_v\tan\theta}{x_v+y_v\tan\theta}$ instead of $y=\beta$ and $y=\tan\phi=\frac{y_v}{x_v}$ instead of $y=\phi$. Obviously, the result of the observability analysis does not change.

\vskip.2cm
\subsubsection{First Case: $y=r^2$} 
We apply the analytic criterion in chapter \ref{ChapterEORC}. We obtain:

\subsubsection{Step 1}

We obviously obtain the same expression as in the case $y=r^2$ of section \ref{SubSectionObservabilityKnown}, i.e., $L^1_g=2(x_v\cos\gamma + y_v\sin\gamma)$, which does not vanish, in general.

\subsubsection{Step 2}

We have  $\Omega_0=span\{[x_v,y_v,0,0]\}$ and $\Omega_1=span\{[x_v,y_v,0,0], [\cos\theta, \sin\theta, y_v\cos\theta-x_v\sin\theta,0]~\}$.

\subsubsection{Step 3}

We have $\mathcal{L}_{^1\phi_0}L^1_g=\mathcal{L}_{f^1}L^1_g=2\cos(\gamma-\theta)$, which does not vanish, in general. This suffices to conclude that the considered system is not in the special case considered by lemma \ref{LemmaConvergenceSpecial} and we need to continue with step 4.

\subsubsection{Step 4}

We have $\tau= \frac{1}{2(x_v\cos\gamma + y_v\sin\gamma)^2}$, as in the case $y=r^2$ of section \ref{SubSectionObservabilityKnown}.  On the other hand, the differential of $\tau$ also includes the derivative with respect to $\gamma$, namely:

\[
\mathcal{D}\tau=\frac{-1}{(x_v\cos\gamma + y_v\sin\gamma)^3}
\left[\cos\gamma, ~\sin\gamma, ~0, ~y_v\cos\gamma-x_v\sin\gamma \right]
\]

\subsubsection{Step 5}

We need to compute $\Omega_2$ and, in order to do this, we need to compute $^1\phi_1$ and $^2\phi_1$ through algorithm \ref{AlgoPhi1}. We obtain: 
$^1\phi_1=^2\phi_1=\left[\begin{array}{c}
  0 \\
  0  \\
  0 \\
  0 \\
\end{array}
\right]$. By using algorithm \ref{AlgoO1} we compute $\Omega_2$ and we obtain that its dimension is $3$. Additionally,  it is possible to verify that $\mathcal{D}\tau\in \Omega_2$, meaning that $m'=2$.

\subsubsection{Step 6}

By a direct computation, it is possible to check that $\Omega_3=\Omega_2$ meaning that $m^*=2$ and $\Omega^*=\Omega_2$.

\subsubsection{Step 7}

We conclude that the dimension of the observable codistribution is equal to $3(<4)$ and the state is not weakly locally observable. In particular, since the differential of every state component does not belong to $\Omega_2$, we conclude that no state component is observable.

\vskip.2cm
\subsubsection{Second Case: $y=\tan\beta$} 
We apply the analytic criterion in chapter \ref{ChapterEORC}. We obtain:

\subsubsection{Step 1}

We obviously obtain the same expression as in the case $y=\tan\beta$ of section \ref{SubSectionObservabilityKnown}, i.e., $L^1_g=-\frac{y_v \cos\gamma - x_v \sin\gamma}{x_v^2 \cos^2\theta + 2 \sin\theta \cos\theta x_v y_v  - y_v^2 \cos^2\theta + y_v^2}$, which does not vanish, in general.

\subsubsection{Step 2}

We compute  $\Omega_0$ and $\Omega_1$: their dimension are $1$ and $2$, respectively.

\subsubsection{Step 3}

We have $\mathcal{L}_{^2\phi_0}L^1_g=\mathcal{L}_{f^2}L^1_g\neq 0$, in general. This suffices to conclude that the considered system is not in the special case considered by lemma \ref{LemmaConvergenceSpecial} and we need to continue with step 4.

\subsubsection{Step 4}

We have $\tau=\frac{x_v\cos\gamma + y_v\sin\gamma + x_v\cos(\gamma - 2\theta) - y_v\sin(\gamma - 2\theta)}{y_v\cos\gamma - x_v\sin\gamma}$, as in the case $y=\tan\beta$ of section \ref{SubSectionObservabilityKnown}.  On the other hand, the differential of $\tau$ also includes the derivative with respect to $\gamma$.

\subsubsection{Step 5}

By using algorithm \ref{AlgoO1} we compute $\Omega_2$ and we obtain that its dimension is $3$. Additionally,  it is possible to verify that $\mathcal{D}\tau\in \Omega_2$, meaning that $m'=2$.

\subsubsection{Step 6}

By a direct computation, it is possible to check that $\Omega_3=\Omega_2$ meaning that $m^*=2$ and $\Omega^*=\Omega_2$.

\subsubsection{Step 7}

We conclude that the dimension of the observable codistribution is equal to $3(<4)$ and the state is not weakly locally observable. In particular, since the differential of every state component does not belong to $\Omega_2$, we conclude that no state component is observable.

\vskip.2cm
\subsubsection{Third Case: $y=\tan\phi$} 
We apply the analytic criterion in chapter \ref{ChapterEORC}. We obtain:

\subsubsection{Step 1}

We obviously obtain the same expression as in the case $y=\tan\phi$ of section \ref{SubSectionObservabilityKnown}, i.e., $L^1_g=-\frac{y_v\cos\gamma - x_v\sin\gamma}{x_v^2}$, which does not vanish, in general.

\subsubsection{Step 2}

We compute  $\Omega_0$ and $\Omega_1$: their dimension are $1$ and $2$, respectively.

\subsubsection{Step 3}

We have $\mathcal{L}_{^2\phi_0}L^1_g=\mathcal{L}_{f^2}L^1_g\neq 0$, in general. This suffices to conclude that the considered system is not in the special case considered by lemma \ref{LemmaConvergenceSpecial} and we need to continue with step 4.

\subsubsection{Step 4}

We have $\tau= -\frac{2x_v^4\cos\gamma}{x_v^4\sin\gamma - x_v^3y_v\cos\gamma}$, as in the case $y=\tan\phi$ of section \ref{SubSectionObservabilityKnown}.  On the other hand, the differential of $\tau$ also includes the derivative with respect to $\gamma$.

\subsubsection{Step 5}

By using algorithm \ref{AlgoO1} we compute $\Omega_2$ and we obtain that its dimension is $4$. As a result, we do not need to proceed with the remaining steps, since we can directly conclude that the entire state is weakly locally observable.

\begin{table}
\begin{center}
\begin{tabular}{|c|c|c|}
  \hline
  $\gamma$ & Output & State observability\\
  \hline
  \hline
  known & $y=r$ & yes \\
    \hline
  known & $y=\beta$ & yes \\
    \hline
  known & $y=\phi$ & yes \\
  \hline
  \hline
  unknown & $y=r$ & no \\
    \hline
  unknown & $y=\beta$ & no \\
    \hline
  unknown & $y=\phi$ & yes \\
    \hline
\end{tabular}
\end{center}
\caption{Weak local observability of the state in all the considered scenarios}
\label{Table}
\end{table}

\vskip.6cm
\noindent Table \ref{Table} summarizes the results of the observability analysis carried out in this section. We conclude this section by remarking that these results agree with our expectation. 
By using the observability rank condition in \cite{Her77}, we easily obtain that, in absence of the unknown input, the dimension of the observable codistribution is $2$ for the first two observations ($y=r$ and $y=\beta$) and $3$ for the last one ($y=\phi$). In particular, for the first two observations, all the initial states rotated around the vertical axis are indistinguishable. In other words, in these two cases, the system exhibits a continuous symmetry \cite{TRO11}. In presence of the unknown input, when $\gamma$ is known, the aforementioned system invariance is broken and the entire state becomes observable. When $\gamma$ is unknown, the symmetry still remains (and obviously also concerns the new state component $\gamma$). Note also a very important aspect. The presence of an unknown input improves the observability properties of a system (this regards the case when $\gamma$ is known). In particular, if $w=0$ (absence of unknown input), the state becomes unobservable despite the knowledge of the unknown input (we know that it is zero), while, when $w\neq0$, the state is observable even if $w$ is unknown. Note that having an unknown input equal to zero is an event that occurs with zero probability and our theory accounts this fact since it is based on definition \ref{DefindistinguishableStates}. To this regard, note also that the validity of theorem \ref{TheoremSeparation}, which allows us to introduce the algorithms \ref{AlgoPhi1} and \ref{AlgoO1}, holds when the unknown input is different from $0$.

\subsection{Simulations}\label{SubSectionSimulations2D}

The scope of this section is to show that, by using a very simple estimator based on an Extended Kalman Filter, we obtain results which agree with our observability analysis. 

\subsubsection{Simulated trajectories and robot sensors}

\begin{figure}[htbp]
\begin{tabular}{cc}
\includegraphics[width=.5\columnwidth]{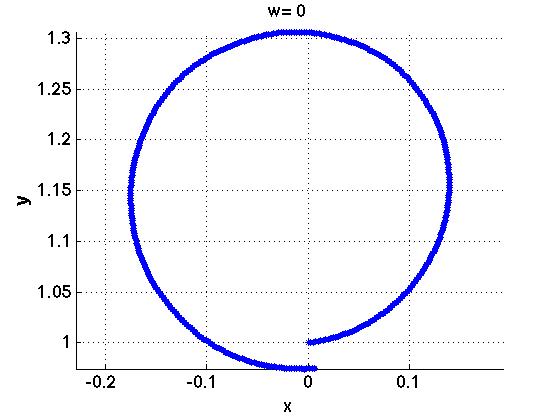}&
\includegraphics[width=.5\columnwidth]{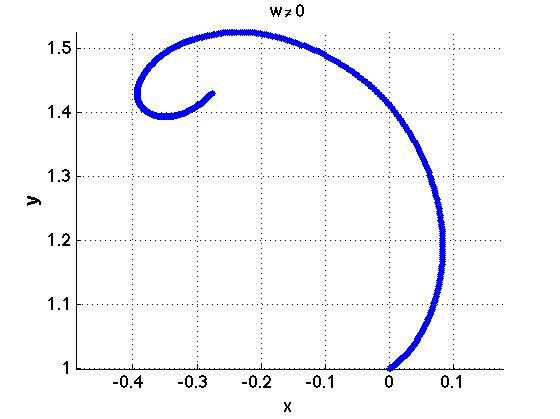}\\
\end{tabular}
\caption{Typical simulated trajectories. Left side without disturbance, right side with disturbance.}
\label{FigTraj1}
\end{figure}

The trajectories are simulated as follows. The equations in (\ref{EquationSystemDynamics}) are discretized with a time step of $5*10^{-4}~s$ and they last for $10~s$. The linear speed, i.e. $v$, is constant and equal to $0.1~m~s^{-1}$. The angular speed, i.e. $\omega$, is generated randomly. Specifically, its value is settled at each step and follows a Gaussian distribution with mean value $\frac{2\pi}{10}~ rad ~s^{-1}$ and variance $\left(\frac{6\pi}{10}\right)^2~ rad^2 ~s^{-2}$. A typical example of trajectory, obtained with this setting, is displayed in fig. \ref{FigTraj1}, left side. The disturbance is generated as follows. The parameter $\gamma$ is set equal to $0.7~\pi$. The unknown input $w$ is generated randomly. Specifically, its value is settled at each step and follows a Gaussian distribution with mean value $0.05 ~m~s^{-1}$ and variance $(0.005)^2~m^2~s^{-2}$. A typical example of trajectory is displayed in fig. \ref{FigTraj1}, right side.

\noindent The robot is equipped with proprioceptive sensors able to measure at each time step the linear and the angular speed. These measurements are affected by errors. Specifically, each measurement is generated at every time step of $5*10^{-4}~s$ by adding to the true value a random error that follows a Gaussian distribution. The mean value of this error is zero and the standard deviation is $0.01$ times the true value for both the linear and the angular speed.

\noindent Regarding the exteroceptive measurements, they are generated at a lower frequency. Specifically, the measurements are generated each $2.5*10^{-2}~s$. Also these measurements are affected by errors. Specifically, the measurement is generated by adding to the true value a random error that follows a Gaussian distribution. The mean value of this error is zero and the standard deviation is $0.01~m$ for the range measurements and $1~deg$ for the two angular measurements (bearing of the origin in the robot frame and bearing of the robot in the global frame).

\subsubsection{Estimation results}

We adopt an Extended Kalman Filter that estimates an extended state that includes the unknown input together with its first order time derivative. In other words, in the case when $\gamma$ is known, the estimated state is: $[x_v,~y_v, ~\theta , ~w, ~w^{(1)}]$. Its dynamics are obtained from (\ref{EquationSystemDynamics}) and are:

\begin{equation}\label{EquationSystemDynamicsExtended}
\left[\begin{aligned}
  \dot{x}_v ~~&= v \cos\theta  +w \cos\gamma\\
  \dot{y}_v ~~&= v \sin\theta  +w \sin\gamma\\
  \dot{\theta} ~~~&= \omega \\
  \dot{w} ~~~&= w^{(1)} \\
  \dot{w}^{(1)} &= w^{(2)} \\
\end{aligned}\right.
\end{equation}

\noindent When $\gamma$ is unknown, the estimated state also includes $\gamma$ and the dynamics are given by (\ref{EquationSystemDynamicsExtended}) and the additional equation $\dot{\gamma}=0$. 

\begin{figure}[htbp]
\begin{tabular}{cc}
\includegraphics[width=.5\columnwidth]{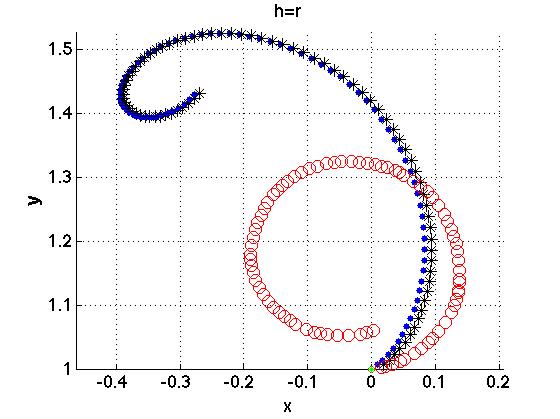}&
\includegraphics[width=.5\columnwidth]{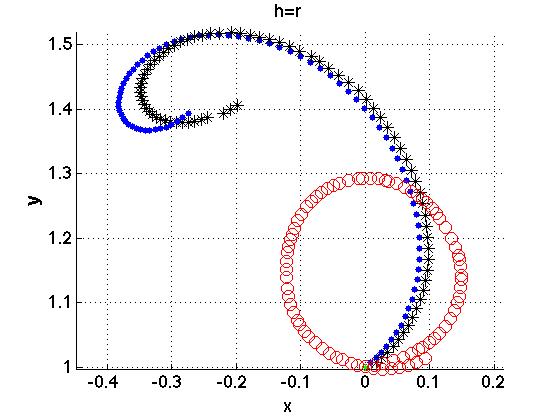}\\
\end{tabular}
\caption{Estimated trajectory in the case when the output is $h=r$. Dots blue are the ground truth, red circles the estimated trajectory by only using the knowledge of the proprioceptive measurements, black stars the trajectory estimated by our estimator. On the left side is the case when $\gamma$ is known and on the right side when it is unknown.}
\label{FigRange}
\end{figure}

\begin{figure}[htbp]
\begin{tabular}{cc}
\includegraphics[width=.5\columnwidth]{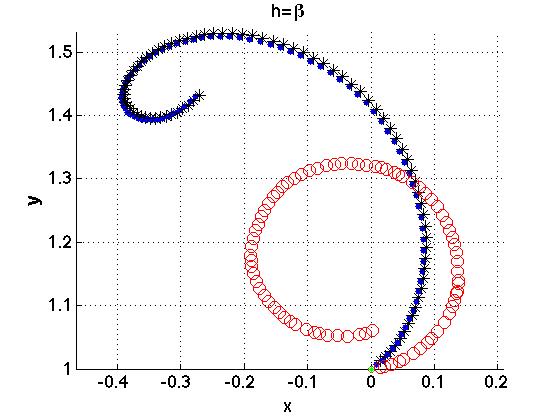}&
\includegraphics[width=.5\columnwidth]{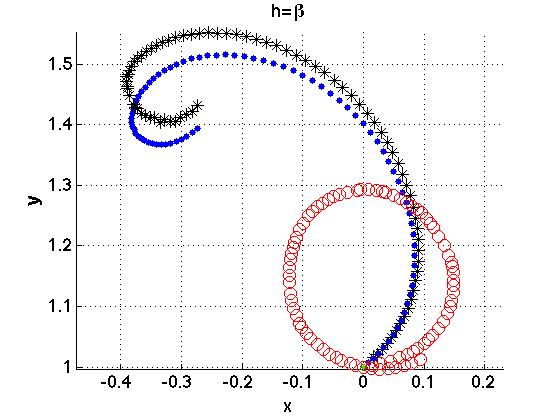}\\
\end{tabular}
\caption{The same as in figure \ref{FigRange} in the case when the output is $h=\beta$.}
\label{FigBeta}
\end{figure}

\begin{figure}[htbp]
\begin{tabular}{cc}
\includegraphics[width=.5\columnwidth]{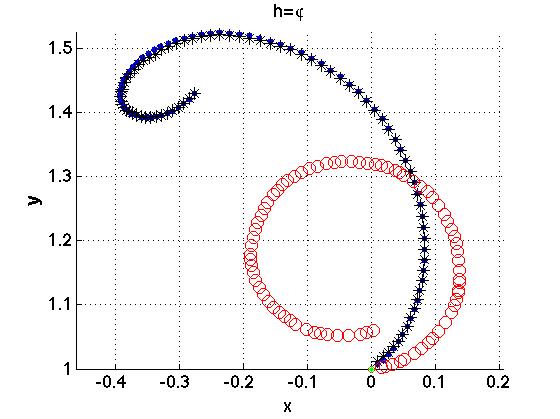}&
\includegraphics[width=.5\columnwidth]{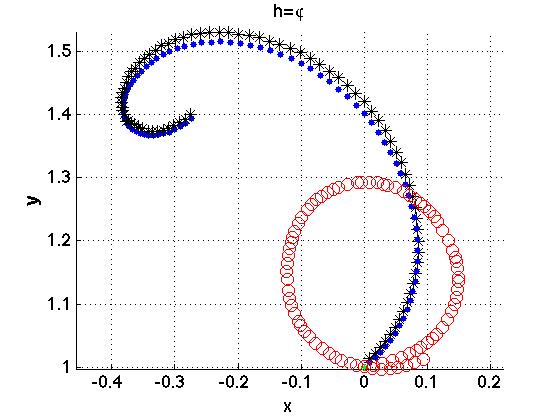}\\
\end{tabular}
\caption{The same as in figure \ref{FigRange} in the case when the output is $h=\phi$.}
\label{FigPhi}
\end{figure}

In order to implement the prediction phase of our Extended Kalman Filter we have to provide the value of $w^{(2)}$, which is unknown. We set this quantity to zero. Note that the simulated trajectory does not satisfy this hypothesis since the disturbance is randomly generated. However, the estimator is able to provide good performance as it is shown in figures \ref{FigRange}-\ref{FigPhi}.

Figure \ref{FigRange} displays the estimated trajectory in the case when the output is provided by the range sensor ($h=r$). Dots blue are the ground truth, red circles the estimated trajectory by only using the knowledge of the proprioceptive measurements (i.e., the measurements of $v$ and $\omega$), black stars the trajectory estimated by our estimator. Left side is the case when $\gamma$ is known and right side when it is unknown. In accordance with our observability analysis, the estimator follows the true trajectory only in the case when $\gamma$ is known.
Figure \ref{FigBeta} displays the estimated trajectory in the case when the output is provided by the on-board bearing sensor ($h=\beta$). Also in this case, the estimator follows the true trajectory only in the case when $\gamma$ is known, in accordance with our observability analysis.
Finally, figure \ref{FigPhi} displays the estimated trajectory in the case when the output is provided by the bearing sensor at the origin ($h=\phi$). In this case, the estimator follows the true trajectory both in the case when $\gamma$ is known and when it is unknown. Also this result agrees with our observability analysis.

\noindent We conclude this section by remarking that the results of our simulations fully agree with the observability analysis provided in sections \ref{SubSectionObservabilityKnown} and \ref{SubSectionObservabilityUnknown}.
Note that we showed the result of a single trial. However, we found always similar results by running many trials.

\section{Vehicle moving in 3D in presence of a disturbance}\label{SectionApplicationDefinition}

\subsection{The system}

\begin{figure}[htbp]
\begin{center}
\includegraphics[width=.8\columnwidth]{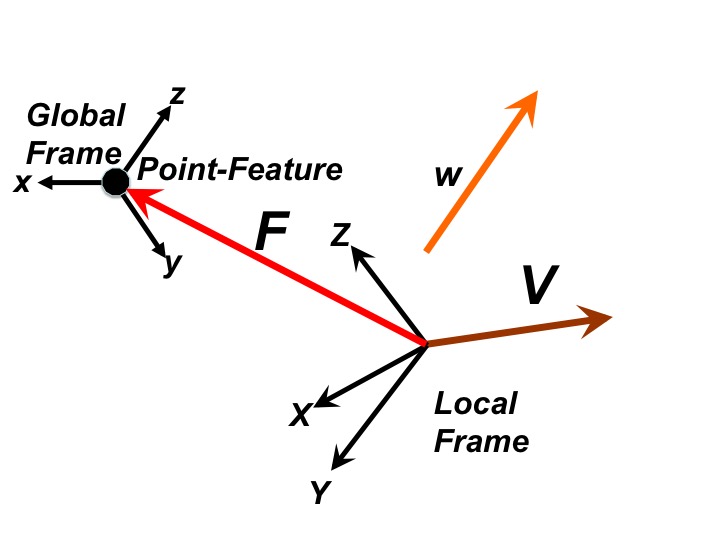}
\caption{Local and global frame for the considered problem. The $z-$axis of the latter is aligned with the direction of the disturbance (assumed to be known and constant in time). The speed $V$ is the vehicle speed with respect to the air, which differs from the ground speed because of the disturbance ($w$).} \label{FigVAir3D}
\end{center}
\end{figure}

We consider a vehicle that moves in a $3D-$environment. We assume that the dynamics of the vehicle are affected by the presence of a disturbance (e.g., this could be an aerial vehicle in presence of wind). We assume that the direction of the disturbance is constant in time and a priori known.  Conversely, the disturbance magnitude is unknown and time dependent. The vehicle is equipped with speed sensors (e.g., airspeed sensors in the case of an aerial vehicle), gyroscopes and a bearing sensor (e.g., monocular camera). We assume that all the sensors share the same frame (in other words, they are extrinsically calibrated). Without loss of generality, we define the vehicle local frame as this common frame. The airspeed sensors measure the vehicle speed with respect to the air in the local frame. The gyroscopes provide the angular speed in the local frame. Finally, the bearing sensor provides the bearing angles of the features in the environment expressed in its own local frame. 
We consider the extreme case of a single point feature and, without loss of generality, we set the origin of the global frame at this point feature (see figure \ref{FigVAir3D} for an illustration). Additionally,  we assume that the $z-$axis of the global frame is aligned with the direction of the disturbance.

\noindent Our system can be characterized by the following state:

\begin{equation}\label{EquationApplicationVAir3DState}
X \triangleq  \left[x, ~y, ~z,  ~q_t, ~q_x, ~q_y, ~q_z \right]^T
\end{equation}

\noindent where $r=[x, ~y, ~z]$ is the position of the vehicle in the global frame and $q=q_t+q_xi+q_yj+q_z k$ is the unit quaternion that describes the transformation change between the global and the local frame. The dynamics are affected by the presence of the disturbance. The disturbance is characterized by the following vector (in the global frame):

\begin{equation}\label{EquationApplicationVAir3DWind}
\bar{w}=w\left[\begin{array}{c}
0\\
0\\
1\\
\end{array}
\right]
\end{equation}

\noindent where $w$ is its unknown magnitude. 

In the sequel, 
for each vector defined in the $3D$ space, the pedix $q$ will be adopted to denote the corresponding imaginary quaternion. 
For instance, regarding the vehicle position, we have: $r_q =0+x~i + y~j + z ~k$.
Additionally, we denote by $V$ and $\Omega$ the following physical quantities:

\begin{itemize}

\item $V=[V_x, ~V_y, ~V_z]$ is the vehicle speed with respect to the air expressed in the local frame (hence, $w~k + q V_q q^*$ is the vehicle speed with respect to the ground expressed in the global frame).

\item $\Omega \triangleq \left[\Omega_x~\Omega_y~\Omega_z \right]$ is the angular speed (and $\Omega_q=0+\Omega_x~i + \Omega_y~j + \Omega_z ~k$).

\end{itemize}

\noindent The dynamics of the state are:

\begin{equation}\label{EquationApplicationVAir3DDynamics}
\left[\begin{aligned}
\dot{r}_q ~~~~~&= w~k + q V_q q^*\\
\dot{q} ~~~~~~&= \frac{1}{2}q\Omega_q\\
\end{aligned}
\right.
\end{equation}

\noindent The monocular camera provides the position of the feature in the local frame ($F_q=-q^*r_qq$) up to a scale. Hence, it provides the ratios of the components of $F$:

\begin{equation}\label{EquationApplicationVAir3DOutputCam}
h_{cam}(X)\triangleq  [h_u,~h_v]^T=\left[
\frac{(q^*r_qq)_x}{(q^*r_qq)_z},
~\frac{(q^*r_qq)_y}{(q^*r_qq)_z} \right]^T
\end{equation}

\noindent where the pedices $x$, $y$ and $z$ indicate
respectively the $i$, $j$ and $k$ component of the corresponding
quaternion. We have also to consider the constraint $q^*q=1$. This provides the further observation:

\begin{equation}\label{EquationApplicationVAir3DOutputq}
h_{const}(X)\triangleq  q^*q
\end{equation}

\noindent Our system is characterized by the state in (\ref{EquationApplicationVAir3DState}), the dynamics in (\ref{EquationApplicationVAir3DDynamics}) and the three outputs $h_u$, $h_v$ and $h_{const}$ in (\ref{EquationApplicationVAir3DOutputCam}) and (\ref{EquationApplicationVAir3DOutputq}).

%
%

\subsection{Observability in absence of disturbance}\label{SubSectionObsNoWind}

Our system is characterized by the state in (\ref{EquationApplicationVAir3DState}), the dynamics in (\ref{EquationApplicationVAir3DDynamics}) with $w=0$ and the three outputs $h_u$, $h_v$ and $h_{const}$ in (\ref{EquationApplicationVAir3DOutputCam}) and (\ref{EquationApplicationVAir3DOutputq}).

By comparing (\ref{EquationApplicationVAir3DDynamics}) with (\ref{EquationSystemDefinitionDynamics1}) we obtain that our system is characterized by
six known inputs ($m_u=6$) that are: $u_1=\Omega_x$, $u_2=\Omega_y$, $u_3=\Omega_z$, $u_4=V_x$ , $u_5=V_y$ and $u_6=V_z$. Additionally, we obtain:

\[
f^1=\frac{1}{2}\left[\begin{array}{c}
 0 \\
 0  \\
 0  \\
 -q_x\\
  q_t\\
  q_z\\
 -q_y\\
\end{array}
\right], ~f^2=\frac{1}{2}\left[\begin{array}{c}
 0 \\
 0  \\
 0  \\
 -q_y\\
 -q_z\\
  q_t\\
  q_x\\
\end{array}
\right], ~f^3=\frac{1}{2}\left[\begin{array}{c}
 0 \\
 0  \\
 0  \\
 -q_z\\
  q_y\\
 -q_x\\
  q_t\\
\end{array}
\right],
\]

\[
f^4=\left[\begin{array}{c}
 q_t^2 + q_x^2 - q_y^2 - q_z^2\\
         2q_tq_z + 2q_xq_y\\
         2q_xq_z - 2q_tq_y\\
                         0\\
                         0\\
                         0\\
                         0\\
\end{array}
\right], ~f^5=\left[\begin{array}{c}
         2q_xq_y - 2q_tq_z\\
 q_t^2 - q_x^2 + q_y^2 - q_z^2\\
         2q_tq_x + 2q_yq_z\\
                         0\\
                         0\\
                         0\\
                         0\\
\end{array}
\right],
\]

\[
f^6=\left[\begin{array}{c}
         2q_tq_y + 2q_xq_z\\
         2q_yq_z - 2q_tq_x\\
 q_t^2 - q_x^2 - q_y^2 + q_z^2\\
                         0\\
                         0\\
                         0\\
                         0\\
\end{array}
\right]
\]

\noindent Finally, in absence of disturbance we have:

\[
g=\left[0, ~0, ~0, ~0, ~0, ~0, ~0 \right]^T
\]


\noindent In this case we can apply the observability rank condition, i.e., algorithm \ref{AlgoHK}, to obtain the observable codistribution.
We compute the codistribution $\Omega_0$ by computing the differentials of the three functions $h_u$, $h_v$ and $h_{const}$. We obtain that this codistribution has dimension equal to $3$.
We use algorithm \ref{AlgoHK} to compute $\Omega_1$. We obtain that its dimension is $4$. In particular, the additional covector is obtained by the differential of the following Lie derivative:

\[
\mathcal{L}_{f^4} h_u
\]

\noindent In other words:

\[
\Omega_1= span\left\{\mathcal{D} h_u, ~\mathcal{D} h_v, ~\mathcal{D} h_{const}, ~\mathcal{D}  \mathcal{L}_{f^4}h_u,\right\}
\]

\noindent All the remaining first order Lie derivatives have differential that is in the above codistribution. Additionally, by an explicit computation, it is easy to realize that 
$\Omega_2=\Omega_1$. This means that algorithm \ref{AlgoHK} has converged and the observable codistribution is $\Omega_1$.

%
%
%
%
%
%
%
%
%
%
%
%
%
%
%

\noindent  By an explicit computation, it is possible to check that the differentials of the components of the vector $F$ belong to $\Omega_1$. This means that all the observable modes are the components of $F$, i.e., the position of the feature in the local frame (obviously, the fourth observable mode is the norm of the quaternion). In particular, no component of the vehicle orientation is observable.

\subsection{Observability in presence of the disturbance}\label{SubSectionObsWind}

We now consider the case when the dynamics are affected by the presence of the disturbance. 
By comparing (\ref{EquationApplicationVAir3DDynamics}) with (\ref{EquationSystemDefinitionDynamics1}) we obtain that the vector fields that characterize the dynamics are the same that characterize the dynamics in absence of disturbance with the exception of the last one, which becomes:

\[
g=\left[0, ~0, ~1, ~0, ~0, ~0, ~0 \right]^T
\]

%
%
%

\noindent To derive the observability properties we apply the proposed analytic tool, by following the seven steps provided in chapter \ref{ChapterEORC}.

\vskip.2cm
\noindent {\bf First Step}

We start by computing the Lie derivatives of the outputs $h_u$, $h_v$ and $h_{const}$ along the vector field $g$. We find that the result differs from zero for the first two outputs. Hence, we use the first output ($h_u$) to define $L^1_g$ (we could choose also the second output $h_v$). In particular, we obtain: $L^1_g \triangleq  \mathcal{L}_g h_u=$

\[
\frac{-y(2q_tq_z - 2q_xq_y) - x(q_t^2 - q_x^2 + q_y^2 - q_z^2)}{[z(q_t^2- q_x^2- q_y^2+ q_z^2)  + 2x(q_tq_y + q_xq_z)  + 2y(q_yq_z  - q_tq_x)]^2}
\]

\vskip.2cm
\noindent {\bf Second Step}

We compute the codistribution $\Omega_0$ by computing the differentials of the three functions $h_u$, $h_v$ and $h_{const}$. This coincides with the case without disturbance, and we obtain that this codistribution has dimension equal to $3$.

We use algorithm \ref{AlgoO1} to compute $\Omega_1$. We obtain that its dimension is $5$. In particular, the additional two independent covectors are obtained by the differentials of the following two Lie derivatives:

\[
\mathcal{L}_{f^4} h_u, ~~~~\mathcal{L}_{\frac{g}{L^1_g}} h_v
\]

\noindent In other words:

\[
\Omega_1= span\left\{\mathcal{D} h_u, ~\mathcal{D} h_v, ~\mathcal{D} h_{const}, ~\mathcal{D}  \mathcal{L}_{f^4}h_u, ~\mathcal{D} \mathcal{L}_{\frac{g}{L^1_g}} h_v\right\}
\]

\noindent All the remaining first order Lie derivatives have differential that is in the above codistribution.

\vskip.2cm
\noindent {\bf Third Step}

We compute $^1\phi_1$, $^2\phi_1$, $^3\phi_1$, $^4\phi_1$, $^5\phi_1$ and $^6\phi_1$ by using algorithm \ref{AlgoPhi1}. We obtain that all these vectors vanish. As a result, all the subsequent steps of algorithm \ref{AlgoPhi1} provide null vectors. Therefore, the assumptions of lemma \ref{LemmaConvergenceSpecial} are trivially met. We set $m'=0$ and we skip to the sixth step.

\vskip.2cm
\noindent {\bf Sixth Step}

We use algorithm \ref{AlgoO1} to compute $\Omega_2$ and we obtain:

\[
\Omega_2=\Omega_1+ span\left\{
~\mathcal{D} \mathcal{L}_{f^4}\mathcal{L}_{\frac{g}{L^1_g}} h_v
\right\}
\]

\noindent Hence, its dimension is $6$. Finally, by using again algorithm \ref{AlgoO1} it is possible to compute $\Omega_3$ and to check that $\Omega_3=\Omega_2$. This means that the algorithm has converged and the observable codistribution is $\Omega^*=\Omega_2$.

\vskip.2cm
\noindent {\bf Seventh Step}

By computing the distribution orthogonal to the codistribution $\Omega^*$ we can find the continuous symmetry that characterizes the unobservable space \cite{TRO11}. By an explicit computation we obtain the following vector: 

\[
\left[ -y, ~ x,~ 0,~ -\frac{q_z}{2},~ -\frac{q_y}{2},~ \frac{q_x}{2},  \frac{q_t}{2} \right]^T
\]


\noindent  This symmetry corresponds to an invariance with respect to a rotation around the $z-$axis of the global frame. This means that we have a single unobservable mode that is the yaw in the global frame\footnote{Note that the chosen global frame is aligned with he direction of the disturbance (fig. \ref{FigVAir3D}). Hence, what is unobservable is a rotation around the direction of the disturbance.}. We conclude by remarking that the presence of the disturbance, even if its magnitude is unknown and is not constant, makes observable the roll and the pitch angles. This result is similar to the result that we obtain in the case of visual and inertial sensor fusion in presence of gravity. The presence of gravity makes observable the roll and the pitch angles, even if its magnitude is unknown \cite{TRO12,FnT14}. What it is non intuitive in the case now investigated, is that, not only the magnitude of the disturbance is unknown, but it is also time dependent.

\subsection{Simulations}\label{SubSectionSimulations3D}

As in section \ref{SubSectionSimulations2D}, the scope of this section is to show that, by using a simple estimator based on an Extended Kalman Filter, we obtain results that agree with our observability analysis.

\vskip.2cm

\noindent {\bf Simulated trajectories and robot sensors}

The trajectories are simulated as follows. The equations in (\ref{EquationApplicationVAir3DDynamics}) are discretized with a time step of $0.01~s$. Each trial lasts $200~s$. The initial vehicle speed is set to zero. The initial position is set equal to $[0.5 ~0.5 ~0.5] m$. 
The vehicle motion is randomly generated. 
The angular speed, i.e. $\Omega$, is Gaussian. Specifically, its value at each step follows a zero mean Gaussian distribution with covariance matrix equal to $(1~deg)^2 I_3$, where $I_3$ is the identity $3\times 3$ matrix.
At each time step, the vehicle speed is incremented by adding a random vector with zero mean Gaussian distribution. In particular, the covariance matrix of this distribution is set equal to $\sigma^2 I_3$, with $\sigma=10^{-4}m$. Finally, we simulate the disturbance as a vector along the $z-$axis, whose magnitude is generated as a random Gaussian variable, with mean value $0.2 ~ms^{-1}$ and variance $0.05^2~m^2s^{-2}$.

\begin{figure}[htbp]
\includegraphics[width=.9\columnwidth]{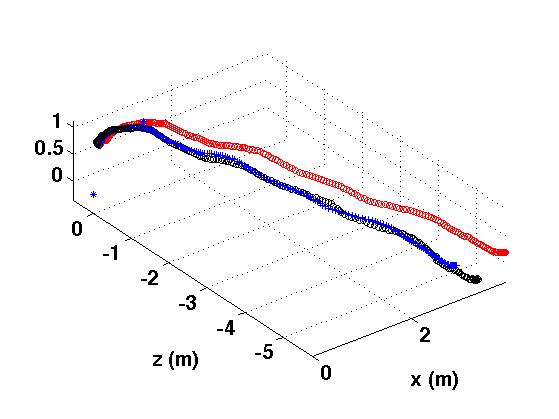}
\caption{A typical simulated trajectory in absence of disturbance. In blue the true trajectory, in red the trajectory estimated by only integrating the gyroscope and the airspeed measurements and in black the trajectory estimated by the proposed EKF.}
\label{FigTrajNo}
\end{figure}

\begin{figure}[ht]\centering
\begin{tikzpicture}[
    zoomboxarray
]
    \node [image node] { \includegraphics[width=0.3\textwidth]{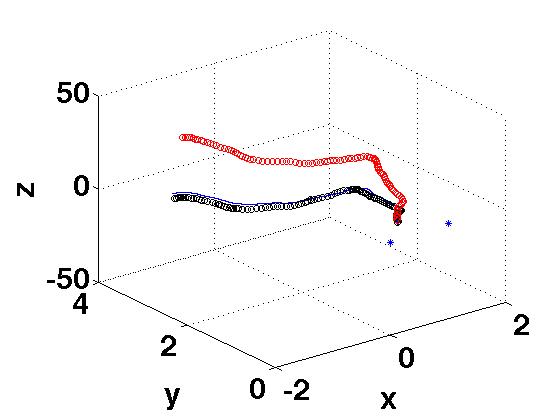} };
    \zoombox{0.32,.53}
\end{tikzpicture}
\caption{A typical simulated trajectory in presence of disturbance along the $z-$axis. In blue the true trajectory, in red the trajectory estimated by only integrating the gyroscope and the airspeed measurements and in black the trajectory estimated by the proposed EKF.}
\label{FigTraj}
\end{figure}

Typical trajectories, obtained with this setting, are displayed in figures \ref{FigTrajNo} and \ref{FigTraj}.

\noindent The vehicle is equipped with proprioceptive sensors able to measure at each time step the speed with respect to the air and the angular speed. These measurements are affected by errors. Specifically, each measurement is generated at every time step of $0.01~s$ by adding to the true value a random error that follows a Gaussian distribution. The mean value of this error is zero and the standard deviation is $0.01~ms^{-1}$ for the airspeed and $1~deg~ s^{-1}$ for the gyroscope.

Regarding the exteroceptive measurements, they are generated at a lower frequency. Specifically, the measurements are generated each $0.1~s$. Also these measurements are affected by errors. Specifically, each measurement is generated by adding to the true value a random error that follows a zero mean Gaussian distribution, with variance $1~deg^2$. We simulate two point features, one at the origin and the other one at the position $[1 ~0 ~1] m$ (note that the result of the observability analysis is independent of the number of point features).

\vskip.2cm

\noindent {\bf Estimation results}

We adopt an EKF that estimates an extended state that includes the state in (\ref{EquationApplicationVAir3DState}) and the unknown input ($w$) together with its first order time derivative (${w}^{(1)}\triangleq \frac{dw}{dt}$). In other words, the estimated state is: $\left[x, ~y, ~z,  ~q_t, ~q_x, ~q_y, ~q_z, ~w, ~{w}^{(1)}  \right]^T$. Its dynamics are given in (\ref{EquationApplicationVAir3DDynamics}) with the two further equations:

\begin{equation}\label{EquationDynamicsAdd}
\left[\begin{aligned}
  \dot{w} ~~~&= w^{(1)} \\
  \dot{w}^{(1)} &= w^{(2)} \\
\end{aligned}\right.
\end{equation}

\begin{figure}[htbp]
\includegraphics[width=.8 \columnwidth]{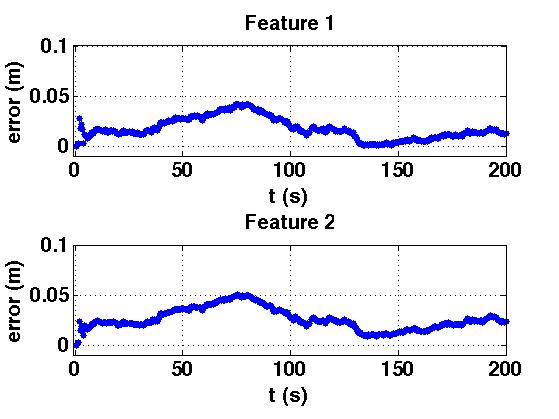}
\caption{Estimation error on the distance of the two observed point features in absence of the disturbance.}
\label{FigFNo}
\end{figure}

\begin{figure}[htbp]
\includegraphics[width=.8 \columnwidth]{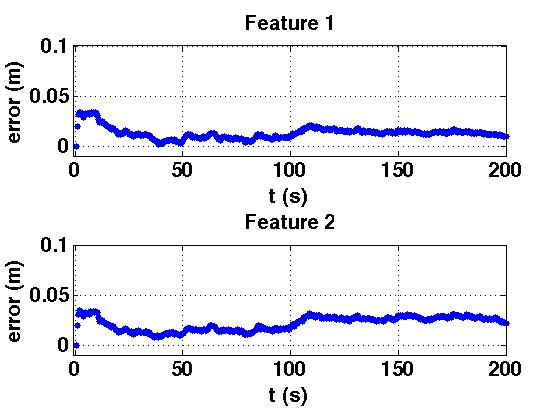}
\caption{As in fig. \ref{FigFNo} but in presence of the disturbance.}
\label{FigF}
\end{figure}

\begin{figure}[htbp]
\begin{tabular}{cc}
\includegraphics[width=.42\columnwidth]{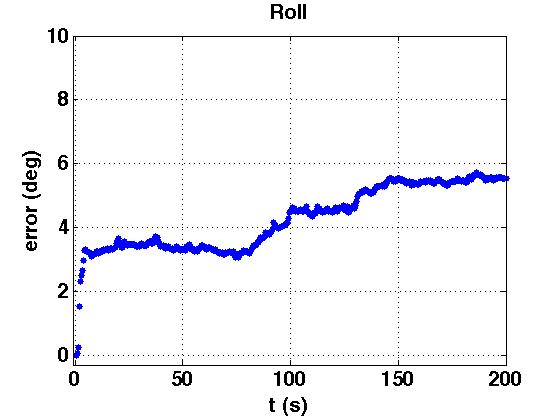}&
\includegraphics[width=.42\columnwidth]{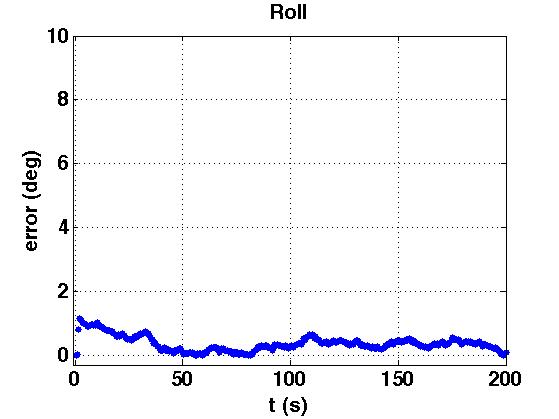}\\
\end{tabular}
\caption{Estimation error on the  vehicle roll angle. Left side in absence of the disturbance, right side in presence of the disturbance.}
\label{FigRoll}
\end{figure}

\begin{figure}[htbp]
\begin{tabular}{cc}
\includegraphics[width=.42\columnwidth]{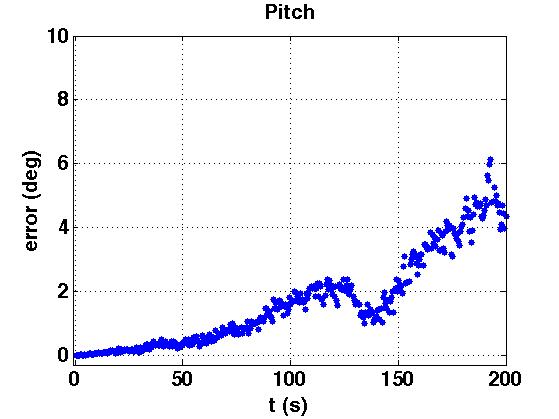}&
\includegraphics[width=.42\columnwidth]{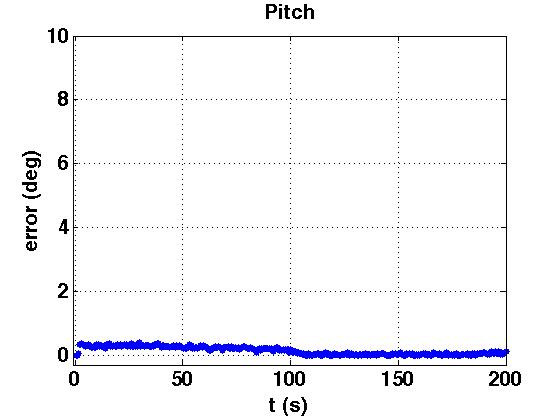}\\
\end{tabular}
\caption{As in fig. \ref{FigRoll} but for the pitch angle.}
\label{FigPitch}
\end{figure}

\begin{figure}[htbp]
\begin{tabular}{cc}
\includegraphics[width=.4\columnwidth]{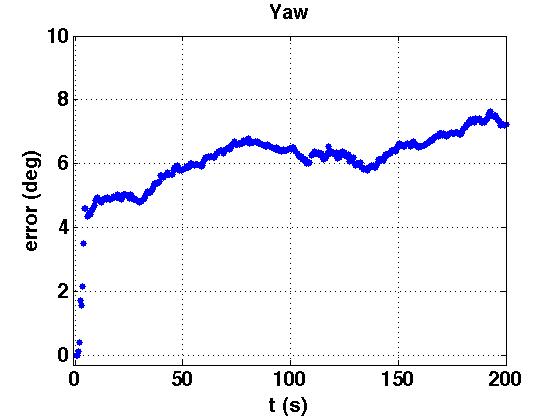}&
\includegraphics[width=.4\columnwidth]{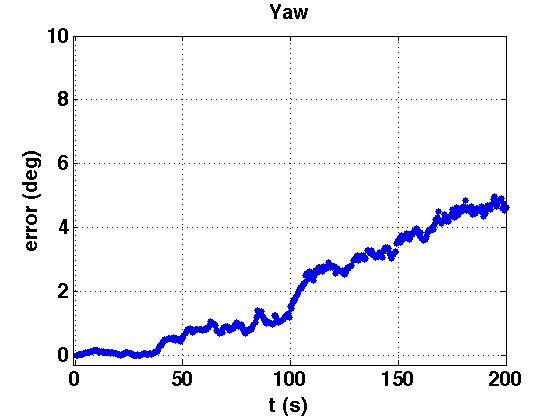}\\
\end{tabular}
\caption{As in fig. \ref{FigRoll} but for the yaw angle.}
\label{FigYaw}
\end{figure}

\noindent In order to implement the prediction phase of our EKF we have to provide the value of $w^{(2)}$ ($\triangleq \frac{d^2w}{dt^2}$), which is unknown. We set this quantity to zero. Note that the simulated trajectory does not satisfy this hypothesis since the disturbance is randomly generated. However, the estimator is able to provide good performance as it is shown in figures \ref{FigFNo}-\ref{FigYaw}.

Figure \ref{FigTrajNo} displays the estimated trajectory in absence of the disturbance. The red line is the trajectory obtained by only using the proprioceptive measurements, while the black line is the trajectory obtained by the proposed EKF. The true trajectory is in blue. Both the estimated trajectories diverge. However, the divergence is slower for the trajectory estimated by the EKF. This divergence is consistent with our observability analysis. The vehicle position is not observable (both with and without disturbance).
Figure \ref{FigTraj} displays the same of figure \ref{FigTrajNo} for a trial in presence of the disturbance. Note that, in this case, the $z-$axis has now a different scale (much larger). This is to show the trajectory obtained by only using the proprioceptive measurements, whose divergence is much faster due to the presence of the disturbance along the $z-$axis.

In the remaining figures, we show the performance of the proposed EKF in estimating the distances of the two features and the roll, the pitch and the yaw angles. In accordance with our observability analysis, the distances of the features are observable both with and without the disturbance. In contrast, the roll and the pitch angles are observable only in presence of the disturbance. Finally, the yaw angle is always unobservable.
All these results are clearly confirmed by our simulation results.

\noindent Figures \ref{FigFNo} and \ref{FigF} show the error on the two distances, in absence and in presence of the disturbance, respectively. We note that, the error does not increase with time in both cases. 

Figures \ref{FigRoll}, \ref{FigPitch} and \ref{FigYaw} show the error on the roll, pitch and yaw angles, respectively. The left side is in absence of the disturbance while the right side is in presence of the disturbance. In absence of the disturbance, the error increases with time for all the three angles. In presence of the disturbance, it does not increase for the roll and the pitch angle.

\noindent We conclude this section by remarking that the results of our simulations fully agree with the observability analysis provided in section \ref{SubSectionObsWind}. 
Also note that we showed the result of a single trial. However, we found always similar results by running many trials.

\chapter{Analytic Derivations}\label{ChapterProofs}

In this chapter, we prove the validity of the analytic results presented in chapter \ref{ChapterObsCod}. We start by introducing several important concepts starting by providing the definition of 
state indistinguishability in presence of unknown inputs (section \ref{SectionProofInd}). These concepts hold for any number of unknown inputs and when the dynamics are affine in the inputs (and not simply linear). Then, section \ref{SectionProof} 
is devoted to prove the results that hold in the case of a single unknown input and therefore to prove the validity of the results provided in chapters \ref{ChapterObsCod} and \ref{ChapterEORC}.

\newpage

\section{State Indistinguishability, State augmentation and Observable Codistribution}\label{SectionProofInd}

As mentioned in the introduction of this chapter, the concepts introduced in this section, hold
for any number of unknown inputs ($m_w$) and when the dynamics are affine in the inputs. Specifically, we will refer to the following dynamics:

\begin{equation}\label{EquationSystemDefinitionDynamicsG}
\begin{aligned}
  \dot{x} &=   g^0(x)+\sum_{i=1}^{m_u}f^i ( x ) u_i +  \sum_{j=1}^{m_w}g^j ( x ) w_j  \\
\end{aligned}
\end{equation}

\noindent where  $g^0( x )$, $ f^i ( x )$, $i=1,\cdots,m_u$, and $ g^j ( x )$, $j=1,\cdots,m_w$, are vector fields in $M\subset \mathbb{R}^n$.

In \cite{Her77,Isi95} the observability properties of a nonlinear system driven by only known inputs are derived starting from the definition of indistinguishable states. According to this definition, it is proven the following fundamental property:

\begin{pr}\label{PropositionHK}
The Lie derivatives of any output computed along any direction allowed by the system dynamics take the same values at the states which are indistinguishable. 
\end{pr}

Starting from this fundamental property it is possible to prove that algorithm \ref{AlgoHK} generates the observable codistribution \cite{Isi95}.
In presence of unknown inputs, we first need to introduce a new definition of indistinguishability. 
To be conservative, the new definition must consider two states indistinguishable if there exists at least one pair of unknown inputs, such that, the outputs obtained starting from the first state under the effect of the first unknown input on a given time interval ($\mathcal{I}$), coincide with the outputs obtained starting from the second state under the effect of the second unknown input on the same time interval, and this holds for any choice of the known inputs. On the other hand, if this condition is achieved only for a unique pair of unknown inputs, the probability that this event occurs is zero. For this reason, we require that this condition is met for infinite pairs of unknown inputs. The concept of measure in spaces with infinite dimensions is not trivial. 
In particular, there is no analogue of Lebesgue measure on an infinite-dimensional Banach space. One possibility, which is frequently adopted, is to use the concept of prevalent and shy sets \cite{Hunt92}. Let us denote by $\mathcal{W}$ the functional space of all the possible unknown input functions. The probability that a given unknown input belongs to a shy subset of $\mathcal{W}$, is $0$. The probability that a given unknown input belongs to a prevalent subset of $\mathcal{W}$, is $1$. Finally, the probability that a given unknown input belongs to a non-shy subset of $\mathcal{W}$, is strictly larger than $0$. 


We introduce the following definition:

\begin{df}[Indistinguishable states]\label{DefindistinguishableStates}
Two states $x_a$ and $x_b$ are indistinguishable if, for every $u(t)$ (the known input vector function), there exists a non-shy subset $\mathcal{W}_a$  in the functional space of all the possible unknown input functions, such that, for any unknown input function $w_a \in \mathcal{W}_a$ it exists an unknown input function $w_b$ such that
$h( x(t; ~ x_a; ~u; ~w_a ))=h( x(t; ~ x_b; ~u; ~w_b ))$ $\forall t\in\mathcal{I}$.
\end{df}


The property stated by proposition \ref{PropositionHK} does not hold in presence of unknown inputs.
Our first objective is to extend the original state in order to obtain a similar property for the resulting extended system. This new property will be the one stated by proposition \ref{PrConstantLieDer}.

To obtain such a result, we proceed as follows. We extend the original state by including the unknown inputs together with their time derivatives. Specifically, we denote by $ ^kx $ the extended state that includes the unknown inputs and their time derivatives up to the $(k-1)-$order:

\begin{equation}\label{EquationExtendedState}
 ^kx  \triangleq [ x^T, ~ w^T, ~ w^{(1)~T}, ~\cdots,~ w^{(k-1)~T}]^T
\end{equation}

\noindent where $w^{(k)} \triangleq \frac{d^k w}{dt^k}$ and $^kx\in M^{(k)}$, with $M^{(k)}$ an open set of $\mathbb{R}^{n+k m_w}$. From (\ref{EquationSystemDefinitionDynamicsG}) it is immediate to obtain the dynamics for the extended state:


\begin{equation}\label{EquationExtendedStateEvolution}
\begin{aligned}
   ^k\dot{x}  &=  G + \sum_{i=1}^{m_u}  F^i u_i + \sum_{j=1}^{m_w}  W^j  w_j^{(k)} \\
\end{aligned}
\end{equation}

\noindent where:

\begin{equation}\label{EquationF_0^E}
G \triangleq \left[\begin{array}{c}
   g^0 + \sum_{j=1}^{m_w}  g^j w_j\\
   w^{(1)} \\
   w^{(2)} \\
  \cdots\\
   w^{(k-1)} \\
   0_{m_w} \\
\end{array}
\right]
\end{equation}

\begin{equation}\label{EquationF_i^E}
F^i \triangleq \left[\begin{array}{c}
   f^i \\
   0_{k m_w} \\
\end{array}
\right], ~~~W^j \triangleq \left[\begin{array}{c}
   0_{n+(k-1) m_w +j-1} \\
   1\\
   0_{m_w-j} \\   
\end{array}
\right]
\end{equation}

\noindent and we denoted by $0_{m}$ the $m-$dimensional zero column vector. We remark that the resulting system has still $m_u$ known inputs and $m_w$ unknown inputs. However, while the $m_u$ known inputs coincide with the original ones,  the $m_w$ unknown inputs are now the $k-$order time derivatives of the original unknown inputs. The state evolution depends on the known inputs via the vector fields $F^i$, ($i=1,\cdots,m_u$) and it depends on the unknown inputs via the unit vectors $W^j$, ($j=1,\cdots,m_w$). Finally, we remark that only the vector field $G$ depends on the new state elements.

In the sequel, we will denote the extended system by \uiok. Additionally, we use the notation: $ \xi  \triangleq [ w^T, ~ w^{(1)~T}, ~\cdots,~ w^{(k-1)~T}]^T$. In this notation we have $ ^kx  = [ x^T, ~ \xi^T]^T$. 
We also denote by \uio $~$ the original system, i.e., the one characterized by the state $x$ and the equations in (\ref{EquationSystemDefinitionDynamicsG}).

We start by providing a simple result for \uiok:

\begin{lm}\label{LmLieDer1}
In \uiok, the Lie derivatives of the output up to the $m^{th}$ order ($m\le k$) are independent of $w_j^{(f)}$, $j=1,\cdots,m_w$, $\forall f\ge m$.
\end{lm}

\proof{ We proceed by induction on $m$ for any $k$. When $m=0$ we only have one zero-order Lie derivative (i.e., $h(x)$), which only depends on $x$, namely it is independent of $w^{(f)}$, $\forall f \ge 0$. Let us assume that the previous assert is true for $m$ and let us prove that it holds for $m+1$. If it is true for $m$, any Lie derivative up to the $m^{th}$ order is independent of $w^{(f)}$, for any $f\ge m$. In other words, the analytical expression of any Lie derivative up to the $m-$order is represented by a function $g(x,w,w^{(1)},\cdots,w^{(m-1)})$. Hence, $d  g = [\frac{\partial g}{\partial x}, \frac{\partial g}{\partial w}, \frac{\partial g}{\partial w^{(1)}}, \cdots, \frac{\partial g}{\partial w^{(m-1)}}, 0_{(k-m)m_w}]$. It is immediate to realize that the product of this differential by any vector field in (\ref{EquationExtendedStateEvolution}) depends at most on $w^{(m)}$, i.e., it is independent of $w^{(f)}$, $\forall f \ge m+1$ $\blacksquare$}

\vskip .2cm

\noindent A simple consequence of this lemma is the following property:

\begin{pr}\label{PrLieDer1}
Let us consider the system \uiok. The Lie derivatives of the output up to the $k^{th}$ order along at least one vector among $W^j$ ($ j=1, \cdots,m_w$) are identically zero.
\end{pr}

\proof{  From the previous lemma it follows that all the Lie derivatives, up to the $(k-1)-$order are independent of $w^{(k-1)}$, which are the last $m_w$ components of the extended state in (\ref{EquationExtendedState}). Then, the proof follows from the fact that any vector among $W^j$ ($ j=1, \cdots,m_w$) has the first $n+(k-1)m_w$ components equal to zero $\blacksquare$}

We have the following property:

\begin{pr}\label{PrLieDer2}
The Lie derivatives of the output up to the $k^{th}$ order along any vector field $G$, $F^1 , \cdots,  F^{m_u}$ for the system \uiok coincide with the same Lie derivatives for the system \uioK
\end{pr}

\proof{ We proceed by induction on $m$ for any $k$. When $m=0$ we only have one zero-order Lie derivative (i.e., $h(x)$), which is obviously the same for the two systems, \uiok and \uioK. Let us assume that the previous assert is true for $m$ and let us prove that it holds for $m+1\le k$. If it is true for $m$, any Lie derivative up to the $m^{th}$ order is the same for the two systems. Additionally, from lemma \ref{LmLieDer1}, we know that these Lie derivatives are independent of $w^{(f)}$, $\forall f \ge m$. The proof follows from the fact that the first $n+mm_w$ components of the vector fields $G$, $F^1 , \cdots,  F^{m_u}$ for \uiok, coincide with the first $n+mm_w$ components of the same vector fields for \uioK, when $m< k$  $\blacksquare$}

\vskip.2cm

\noindent For \uiok $~$ we have a fundamental property that is the extension of the one stated by proposition \ref{PropositionHK}:

\begin{pr}\label{PrConstantLieDer}
If $x_a$ and $x_b$ are indistinguishable, there exist $\xi_a$ and $\xi_b$ such that, in \uiok, the Lie derivatives of the output  up to the $k^{th}$-order, along all the vector fields that characterize the dynamics of \uiok, take the same values at $[x_a, \xi_a]$ and $[x_b, \xi_b]$.
\end{pr}

\proof{We consider a piecewise-constant input $\tilde{u}$ as follows ($i=1,\cdots,m_u$):

\begin{equation}\label{EquationInputPiecewise}
 \tilde{u}_i(t)=
 \end{equation}
\[
 \left\{\begin{aligned}
 &u_i^1 ~~~t\in[0, ~t_1)\\
 &u_i^2 ~~~t\in[t_1, ~t_1+t_2)\\
 &\cdots \\
 &u_i^g ~~~t\in[t_1+t_2+\cdots+t_{g-1}, ~t_1+t_2+\cdots+t_{g-1}+t_g)\\
\end{aligned}\right.
\]

\noindent Since $x_a$ and $x_b$ are indistinguishable, there exist two unknown input functions $w_a(t)$ and $w_b (t)$ such that the output coincide on $x_a$ and $x_b$. In particular, we can write:

\begin{equation}\label{Equationha=hb}
h( x(t; ~ [x_a, \xi_a]; ~\tilde{u};  ~w_a^{(k)} ))=h( x(t; ~ [x_b, \xi_b]; ~\tilde{u};  ~w_b^{(k)} ))
\end{equation}

\noindent $\forall t\in [0, ~t_1+t_2+\cdots+t_{g-1}+t_g)\subset\mathcal{I}$. On the other hand, by taking the two quantities in (\ref{Equationha=hb}) at $t=t_1+t_2+\cdots+t_{g-1}+t_g$, we can consider them as functions of the $g$ arguments $t_1,t_2,\cdots,t_g$. Hence, by differentiating with respect to all these variables, we also have:

\begin{equation}\label{EquationhaD=hbD}
\frac{\partial^g h( x(t_1+\cdots+t_g; ~ [x_a, \xi_a]; ~\tilde{u};  ~w_a^{(k)} ))}{\partial t_1\partial t_2\cdots\partial t_g}=
\end{equation}
\[
=\frac{\partial^g h( x(t_1+\cdots+t_g; ~ [x_b, \xi_b]; ~\tilde{u};  ~w_b^{(k)} ))}{\partial t_1\partial t_2\cdots\partial t_g}
\]

\noindent By computing the previous derivatives at $t_1=t_2=\cdots=t_g=0$ and by using proposition \ref{PrLieDer1} we obtain, if $g\le k$:

\begin{equation}\label{EquationLDa=LDb}
\mathcal{L}^g_{\theta_1 \theta_2 \cdots \theta_g} h\left|_{\begin{aligned}
 &x=x_a\\
 &\xi=\xi_a\\
\end{aligned}}\right.=\mathcal{L}^g_{\theta_1 \theta_2 \cdots \theta_g} h\left|_{\begin{aligned}
 &x=x_b\\
 &\xi=\xi_b\\
\end{aligned}}\right.
\end{equation}

\noindent where $\theta_h=G+\sum_{i=1}^{m_u}F^iu_i^h$, $h=1,\cdots,g$. The equality in (\ref{EquationLDa=LDb}) must hold for all possible choices of $u_1^h,\cdots,u_{m_u}^h$. By appropriately selecting these $u_1^h,\cdots,u_{m_u}^h$, we finally obtain:

\begin{equation}
\mathcal{L}^g_{v_1 v_2 \cdots v_g} h\left|_{\begin{aligned}
 &x=x_a\\
 &\xi=\xi_a\\
\end{aligned}}\right.=\mathcal{L}^g_{v_1 v_2 \cdots v_g} h\left|_{\begin{aligned}
 &x=x_b\\
 &\xi=\xi_b\\
\end{aligned}}\right.
\end{equation}

\noindent where $v_1 v_2 \cdots v_g$ are vector fields belonging to the set $\{ G, F^1,\cdots,F^{m_u}\}$  $\blacksquare$}

\vskip.2 cm

%
%
%

\noindent The main difference between propositions \ref{PropositionHK} and \ref{PrConstantLieDer} is that, in the latter, we cannot consider any order Lie derivative, since the order cannot exceed $k$. This will have important consequences, as we will see. Before discussing this point, we remark that in \cite{Her77} it was also defined the concept of $V-$indistinguishable states, with $V$ a subset of the definition set that includes the two considered states. From this definition and the previous proof we can alleviate the assumptions in the previous proposition. Specifically, we have the following:

\begin{Rm}\label{Remark}
The statement of proposition \ref{PrConstantLieDer} also holds if $x_a$ and $x_b$ are $V-$indistinguishable.
\end{Rm}

\noindent Now let us discuss how we can use the result stated by the proposition \ref{PrConstantLieDer} to investigate the observability properties. Thanks to the results stated by propositions \ref{PrLieDer2} and \ref{PrConstantLieDer} we can easily build a codistribution that is observable. 
In the sequel, we will denote by $d $ the differential with respect to the entire extended state. We will call {\it original state} the vector $x$ and, we remind the reader that we denote by $\mathcal{D}$ the differential with respect to the original state. The observable codistribution for \uiok $~$is the span of the differentials of all the Lie derivatives of the output along $G$, $F^1 , \cdots,  F^{m_u}$ up to the $k$-order. Hence, for any $m\le k$, it is obtained recursively by the following algorithm:

\begin{al}\label{AlgoOT} Observable codistribution for the extended state ($m\le k$)
\begin{enumerate}
\item $\bar{\Omega}_0=span\{d  h\}$;
\item $\bar{\Omega}_m=\bar{\Omega}_{m-1}+\mathcal{L}_{G}\bar{\Omega}_{m-1}+\sum_{i=1}^{m_u}\mathcal{L}_{F^i}\bar{\Omega}_{m-1}$
\end{enumerate}
\end{al}

\noindent It is possible to obtain all the observability properties of the original state starting from the codistributions generated by the previous algorithm. In order to show this, we  introduce the concept of {\it Observable Mode}. This concept can be also useful in the case without unknown inputs.

\noindent Given $ x_0 $, we denote by $I_{x_0}$ the set of all the states $x$ such that $x$ and $x_0$ are indistinguishable. According to the theory of observability, a system is observable in $x_0$ if $I_{x_0}=x_0$. We introduce here the following new definition:

\begin{df}[Observable Mode]\label{DefObservableMode}
A scalar function is observable in $x_0$ if it is constant on $I_{x_0}$. Additionally, it is weakly observable in $x_0$ if it exists an open neighbourhood $B_{x_0}$ such that it is constant on $B_{x_0}\cap I_{x_0}$.

\end{df}

\noindent The basic idea behind definition \ref{DefObservableMode} is the following. Let us suppose that the true initial state is $x_0$. According to the definition of indistinguishable set (which is based on definition \ref{DefindistinguishableStates}), the system only contains the information to establish whether the initial state belongs to the indistinguishable set $I_{x_0}$ or not. Hence, if a scalar function takes the same value on this set, we conclude that the system has enough information to know the value of this scalar function at the initial time. We also remark that definition \ref{DefObservableMode} generalizes the definition of observability. Specifically, a system is observable in $x_0$ when all the components of the state $x$ are observable modes in $x_0$. When the scalar function is only weakly observable at $x_0$, we conclude that the system has enough information to know the value of this scalar function at the initial time, provided that we a priori know that the initial state is {\it sufficiently} close to $x_0$.

The following two propositions generalize theorem $3.1$ and $3.11$ in \cite{Her77}.

\begin{pr}\label{Pr3.1}
Given a scalar function $\theta(x)$, if it exists an integer $k$ such that for the extended system \uiok, $d  \theta \in \bar{\Omega}_k$ in a given $x_0$ and for a given extension $\xi_0$, then $\theta$ is weakly observable in $x_0$.
\end{pr}

\proof{$\theta(x)$ can be expressed in terms of the Lie derivatives of the output along the fields that characterize the dynamics of \uiok, up to the $k-$order. We can write  $\theta(x)=\mathcal{G}(\phi_1(x,\xi), \cdots, \phi_L(x,\xi))$, $\forall x\in B_{x_0}$ where $\phi_1(x,\xi), \cdots, \phi_L(x,\xi)$ are $L$ Lie derivatives among the ones of above and $\mathcal{G}$ is a given function.

Let us consider a given $x$ indistinguishable from $x_0$ and that belongs to $B_{x_0}$. In other words, $x \in B_{x_0}\cap I_{x_0}$. From proposition \ref{PrConstantLieDer}  there are two extensions $\xi_a$ and $\xi_b$, such that, the Lie derivatives $\phi_1(x,\xi), \cdots, \phi_L(x,\xi)$ take the same values on the two extended states $[x, \xi_a]$ and $[x_0, \xi_b]$. Therefore, $\forall x\in B_{x_0} \cap I_{x_0}$  we have: $\theta(x)=\mathcal{G}(\phi_1(x,\xi_a), \cdots, \phi_L(x,\xi_a))= \mathcal{G}(\phi_1(x_0,\xi_b), \cdots, \phi_L(x_0,\xi_b))=\theta(x_0)$ and $\theta$ is weakly observable in $x_0$
$\blacksquare$}

\begin{pr}\label{Pr3.11}
If the scalar function $\theta(x)$ is weakly observable in $x_0$, then $\exists k$ and $\xi_0$ such that $d  \theta \in \bar{\Omega}_k$ a.e. on an open neighbourhood $B_{x_0}$.
\end{pr}

\proof{ We proceed by contradiction. $\forall k$ and $\xi$ and for any open ball centered on $x_0$ with radius $r$ ($B_{x_0}^r$), it exists a set $C_{x_0}^r\subseteq B_{x_0}^r$ with measure strictly larger than zero for which $d  \theta \notin \bar{\Omega}_k$. This means that $\forall k$ and $\xi$ it exists a distribution with dimension larger than zero that is orthogonal to the codistribution $\bar{\Omega}_k$ in $C_{x_0}^r$. This means that it exists at least one vector field $s$, of the same dimension of the original state, such that, $\forall k$ and $\xi$, it exists at least one vector field $s_k$ of dimension $km_w$, such that the vector field $[s^T,s_k^T]^T$ is orthogonal to the codistribution $\bar{\Omega}_k$ in $C_{x_0}^r$.
From this, we obtain that it exists $\epsilon>0$ such that the states $x_0$ and $x_0+\epsilon s$ are indistinguishable. Indeed, it is possible to express the $m-$time derivative of the output at the initial time, for a generic integer $m$, in terms of the Lie derivatives of the output along the fields of \uiok, up to the $m-$order, and the known inputs and their time derivatives (see \cite{FnT14}, sect 2.5.1). Since this holds for any $k$, the $m-$order time derivative of the output at the initial time coincides in $x_0$ and $x_0+\epsilon s$. Since this holds for any order $m$, from the Taylor theorem the output coincides on a given time interval and consequently $x_0$ and $x_0+\epsilon s$ are indistinguishable. But this means that also $\theta(x_0)=\theta(x_0+\epsilon s)$. Hence, $d  \theta \cdot s=0$ a.e. on an open neighbourhood $B_{x_0}$ and $d  \theta \in \bar{\Omega}_k$ a.e. on an open neighbourhood $B_{x_0}$
$\blacksquare$}

\vskip.2cm

Propositions \ref{Pr3.1} and \ref{Pr3.11} state that all the observability properties of the original state are contained in the codistributions generated by algorithm \ref{AlgoOT}. From this, we can easily obtain a sufficient condition for the observability of the original state. Indeed, if on a given $x_0$, the differential of a given component of $x$ (the original state) belongs to $\bar{\Omega}_m$ for a given integer $m\le k$, we can conclude that this state component is weakly locally observable (in $x_0$). If this holds for all the state components, we can conclude that the entire original state is weakly locally observable. 
More in general, we can conclude that, a given scalar function of the original state, is weakly observable in a given point $x_0$ if its differential (computed in $x_0$) belongs to $\bar{\Omega}_m$ (computed in $x_0$) for a given integer $m\le k$.
On the other hand, we remark the following two fundamental differences between algorithms \ref{AlgoHK} and \ref{AlgoOT}:

\begin{enumerate}
\item In the latter, since the state augmentation can be continued indefinitely, we do not have convergence;

\item The latter provides a codistribution that describes simultaneously the observability properties of the original state and its extension.
\end{enumerate}

The goal of the next section is to address these fundamental issues in the case of a single unknown input. In particular, we show that it is possible to directly compute the entire observable codistribution of the original system, namely, without the need of extending the state.

\newpage

\section{Proof of the validity of the analytic criterion}\label{SectionProof}

This section is devoted to the case of a single unknown input (i.e., $m_w=1$) and dynamics linear in the inputs (the vector $g^0$ is the null vector). In other words, we are considering the system characterized by 
(\ref{EquationSystemDefinitionDynamics1}).
When $m_w=1$, the extended state that includes the time derivatives of $w$ up to the $(k-1)-$order is:

\begin{equation}\label{EquationExtendedState1}
 ^kx  \triangleq [ x^T, ~ w, ~ w^{(1)}, ~\cdots,~ w^{(k-1)}]^T
\end{equation}

\noindent where $w^{(j)}\triangleq \frac{d^jw}{dt^j}$. The dimension of the extended state is in this case $n+k$. 
For the clarity sake, let us consider the case of a single known input (i.e., we consider the system in (\ref{EquationSystemDefinitionDynamics1}) with $m_u=1$). 
We provide all the analytic results in this simpler case and then we extend them to the case of $m_u>1$.
From (\ref{EquationSystemDefinitionDynamics1}) it is immediate to obtain the dynamics for the extended state:

\begin{equation}\label{EquationExtendedStateEvolution1}
\begin{aligned}
   ^k\dot{x}  &=    G(^kx)+ F( x ) u + W  w^{(k)} \\
\end{aligned}
\end{equation}

\noindent where:

\begin{equation}\label{EquationFG}
 F\triangleq \left[\begin{array}{c}
    f ( x )\\
   0 \\
   0 \\
  \cdots\\
   0 \\
   0 \\
\end{array}
\right]
~~~ G \triangleq \left[\begin{array}{c}
    g ( x ) w\\
   w^{(1)} \\
   w^{(2)} \\
  \cdots\\
   w^{(k-1)} \\
   0 \\
\end{array}
\right]
~~~ W \triangleq \left[\begin{array}{c}
   0\\
   0 \\
   0 \\
  \cdots\\
   0 \\
   1 \\
\end{array}
\right]
\end{equation}

\noindent and we set $f(x) \triangleq f^1(x)$.

Algorithm \ref{AlgoOT} becomes:

\begin{al}\label{AlgoOT_1} Observable codistribution for the extended state in the case $m_u=1$ ($m\le k$)
\begin{enumerate}
\item $\bar{\Omega}_0=span\{d  h\}$;
\item $\bar{\Omega}_m=\bar{\Omega}_{m-1}+\mathcal{L}_{G}\bar{\Omega}_{m-1}+\mathcal{L}_{F}\bar{\Omega}_{m-1}$
\end{enumerate}
\end{al}

\noindent where $G$ is now the simpler vector field given in (\ref{EquationFG}) and $d $ denotes the differential respect to the state in (\ref{EquationExtendedState1}).

We will address the two fundamental issues mentioned at the end of the previous section.
This is obtained by two separates steps. In the first step (sect. \ref{SubSectionSeparation}) we perform a separation on the observable codistribution defined by algorithm \ref{AlgoOT_1}. This codistribution can be split into two codistributions: the former is the codistribution generated by algorithm \ref{AlgoO1}, once embedded in the extended space, and the latter is the codistribution $L^m$ defined in section \ref{SubSectionSeparation} (see theorem \ref{TheoremSeparation}). We prove (Lemma \ref{LemmaFundamental}) that the second codistribution ($L^m$) can be ignored when deriving the observability properties of the original state. In the second step (sect. \ref{SubSectionConvergence}) we prove that algorithm \ref{AlgoO1} converges in at most $n+2$ steps and we provide the convergence criterion. Finally, in section \ref{SubSectionExt}
we extend the results of the previous two theorems to the case of multiple known inputs ($m_u>1$).

\subsection{Separation}\label{SubSectionSeparation}

For each integer $m$, we generate the codistribution $\Omega_m$ by using algorithm \ref{AlgoO1} (note that here we are considering $m_u=1$). By construction, the generators of $\Omega_m$ are the differentials of scalar functions that only depend on the original state ($x$) and not on its extension. In the sequel, we need to embed this codistribution in $\mathbb{R}^{n+k}$. We will denote by $[\Omega_m,0_k]$ the codistribution made by covectors whose first $n$ components are covectors in $\Omega_m$ and the last components are all zero. Additionally, we will denote by $L^m$ the codistribution that is the span of the Lie derivatives of $d  h$ up to the order $m$ along the vector $G$, i.e., $L^m \triangleq span\{\mathcal{L}^1_Gd  h, \mathcal{L}^2_Gd  h, \cdots, \mathcal{L}^m_Gd  h\}$. We finally introduce the following codistribution:

\begin{df}[$\tilde{\Omega}$ codistribution]\label{DefinitionOmegaTilde}
This codistribution is defined as follows: $\tilde{\Omega}_m \triangleq [\Omega_m,0_k]+ L^m$

\end{df}

\noindent The codistribution $\tilde{\Omega}_m$ consists of two parts. Specifically, we can select a basis that consists of exact differentials that are the differentials of functions that only depend on the original state ($x$) and not on its extension (these are the generators of $[\Omega_m,0_k]$) and the differentials $\mathcal{L}^1_Gd  h, \mathcal{L}^2_Gd  h, \cdots, \mathcal{L}^m_Gd  h$. The second set of generators, i.e., the differentials $\mathcal{L}^1_Gd  h, \mathcal{L}^2_Gd  h, \cdots, \mathcal{L}^m_Gd  h$, are $m$ and, with respect to the first set, they are differentials of functions that also depend on the state extension $\xi=[w,~w^{(1)},\cdots, ~w^{(m-1)}]^T$. We have the following result:

\begin{lm}\label{LemmaFundamental}
Let us denote with $x_j$ the $j^{th}$ component of the state ($j=1,\cdots,n$). We have: $\mathcal{D}x_j \in \Omega_m$ if and only if $d  x_j \in \tilde{\Omega}_m$
\end{lm}

\proof{The fact that $\mathcal{D}x_j \in \Omega_m$ implies that $d  x_j \in \tilde{\Omega}_m$ is obvious since $[\Omega_m,0_k] \subseteq\tilde{\Omega}_m$ by definition. Let us prove that also the contrary holds, i.e., that if $d  x_j \in \tilde{\Omega}_m$ then $\mathcal{D}x_j \in \Omega_m$.
Since $d  x_j \in \tilde{\Omega}_m$ we have $d  x_j=\sum_{i=1}^{N_1}c^1_i \omega_i^1+\sum_{i=1}^{N_2}c^2_i \omega_i^2$, where $\omega_1^1, \omega_2^1, \cdots, \omega_{N_1}^1$ are $N_1$ generators of $[\Omega_m,0_k]$, $\omega_1^2, \omega_2^2, \cdots, \omega_{N_2}^2$ are $N_2$ generators of $L^m$ and $c^1_1,\cdots,c^1_{N_1},c^2_1,\cdots,c^2_{N_2}$ are suitable coefficients. We want to prove that $N_2=0$.

We proceed by contradiction. Let us suppose that $N_2\ge 1$.
We remark that the first set of generators have the last $k$ entries equal to zero, as for $d  x_j$. The second set of generators consists of the Lie derivatives of $d  h$ along $G$ up to the $m$ order. Let us select the one that is the highest order Lie derivative and let us denote by $j'$ this highest order. We have $1\le N_2 \le j'\le m$. By a direct computation, it is immediate to realize that this is the only generator that depends on $w^{(j'-1)}$. Specifically, the dependence is linear by the product $L^1_g w^{(j'-1)}$ (we remind the reader that $L^1_g \neq 0$). But this means that $d  x_j$ has the $(n+j')^{th}$ entry equal to $L^1_g\neq 0$ and this is not possible since $d  x_j=[\mathcal{D}x_j,0_k]$ $\blacksquare$}


\noindent A fundamental consequence of this lemma is that, if we are able to prove that $\tilde{\Omega}_m=\bar{\Omega}_m$, the weak local observability of the original state $x$, can be investigated by only considering the codistribution $\Omega_m$.
In the rest of this section we prove this fundamental theorem, stating that $\tilde{\Omega}_m=\bar{\Omega}_m$.

For a given $m\le k$ we define the vector $\Phi_m\in \mathbb{R}^{n+k}$ by the following algorithm:

\begin{enumerate}
\item $\Phi_0=F$;
\item $\Phi_m=[\Phi_{m-1}, ~G]$
\end{enumerate}

\noindent where now the Lie brackets $[\cdot, \cdot]$ are computed with respect to the extended state, i.e.:

\[
[F,~G] \triangleq \frac{\partial G}{\partial ^kx} F(^kx) - \frac{\partial F}{\partial ^kx} G(^kx)
\]

\noindent By a direct computation it is easy to realize that $\Phi_m$ has the last $k$ components identically null. In the sequel, we will denote by $\Breve{\Phi}_m$ the vector in $\mathbb{R}^n$ that contains the first $n$ components of $\Phi_m$. In other words, $\Phi_m\triangleq [\Breve{\Phi}_m^T,0_k^T]^T$. Additionally, we set $\hat{\phi}_m\triangleq \left[\begin{array}{c}
    \phi_m \\
   0_k \\
\end{array}
\right]$ ($\phi_m$ is defined by algorithm \ref{AlgoPhi1}).

We have the following result:

\begin{lm}\label{LemmaLemma1}
$\mathcal{L}_G \bar{\Omega}_m+\mathcal{L}_{\Phi_m}d  h=\mathcal{L}_G \bar{\Omega}_m+\mathcal{L}_F\mathcal{L}^m_Gd  h$
\end{lm}

\proof{We have $\mathcal{L}_F\mathcal{L}^m_Gd  h=\mathcal{L}_G\mathcal{L}_F\mathcal{L}^{m-1}_Gd  h+\mathcal{L}_{\Phi_1}\mathcal{L}^{m-1}_G d  h$.
 
The first term $\mathcal{L}_G\mathcal{L}_F\mathcal{L}^{m-1}_Gd  h \in \mathcal{L}_G \bar{\Omega}_m$. Hence, we need to prove that
$\mathcal{L}_G \bar{\Omega}_m+\mathcal{L}_{\Phi_m}d  h=\mathcal{L}_G \bar{\Omega}_m+\mathcal{L}_{\Phi_1}\mathcal{L}^{m-1}_G d  h$. We repeat the previous procedure $m$ times. Specifically,  we use  the equality $\mathcal{L}_{\Phi_j}\mathcal{L}^{m-j}_Gd  h=\mathcal{L}_G\mathcal{L}_{\Phi_j}\mathcal{L}^{m-j-1}_Gd  h+\mathcal{L}_{\Phi_{j+1}}\mathcal{L}^{m-j-1}_G d  h$, for $j=1,\cdots,m$, and we remove the first term since $\mathcal{L}_G\mathcal{L}_{\Phi_j}\mathcal{L}^{m-j-1}_Gd  h \in \mathcal{L}_G \bar{\Omega}_m$ $\blacksquare$}

\begin{lm}\label{LemmaRisA}
$\Breve{\Phi}_m=\sum_{j=1}^m c^n_j(\mathcal{L}_Gh, \mathcal{L}^2_Gh, \cdots, \mathcal{L}^m_Gh) \phi_j$, i.e., the vector $\Breve{\Phi}_m$ is a linear combination of the vectors $\phi_j$ ($j=1,\cdots,m$), where the coefficients ($c^n_j$) depend on the state only through the functions that generate the codistribution $L^m$
\end{lm}

\proof{ We proceed by induction. By a direct computation it is immediate to obtain: $\Breve{\Phi}_1=\phi_1 \mathcal{L}_Gh$. 

\noindent {\bf Inductive step:} Let us assume that $\Breve{\Phi}_{m-1}=\sum_{j=1}^{m-1} c_j(\mathcal{L}_Gh, \mathcal{L}^2_Gh, \cdots, \mathcal{L}^{m-1}_Gh) \phi_j$. We have:

\[
\Phi_m=[\Phi_{m-1},~G]=\sum_{j=1}^{m-1}\left[c_j \left[\begin{array}{c}
    \phi_j\\
   0_k \\
\end{array}
\right],~G\right]=
\]
\[
\sum_{j=1}^{m-1}c_j\left[\left[\begin{array}{c}
    \phi_j\\
   0_k \\
\end{array}
\right],~G\right] -
\sum_{j=1}^{m-1}\mathcal{L}_G c_j \left[\begin{array}{c}
    \phi_j\\
   0_k \\
\end{array}
\right]
\]

 We directly compute the Lie bracket in the sum (note that $\phi_j$ is independent of the unknown input $w$ and its time derivatives):

\[
\left[\left[\begin{array}{c}
    \phi_j\\
   0_k \\
\end{array}
\right],~G\right]= \left[\begin{array}{c}
    [\phi_j,~g] w\\
   0_k \\
\end{array}
\right]=\left[\begin{array}{c}
    \phi_{j+1} \mathcal{L}^1_Gh\\
   0_k \\
\end{array}
\right]
\]

Regarding the second term, we remark that $\mathcal{L}_G c_j = \sum_{i=1}^{m-1} \frac{\partial c_j}{\partial (\mathcal{L}_G^ih)} \mathcal{L}_G^{i+1}h$. By setting $\tilde{c}_j=c_{j-1}\mathcal{L}^1_Gh$ for $j=2,\cdots,m$ and $\tilde{c}_1=0$, and by setting $\bar{c}_j=-\sum_{i=1}^{m-1} \frac{\partial c_j}{\partial (\mathcal{L}_G^ih)} \mathcal{L}_G^{i+1}h$ for $j=1,\cdots,m-1$ and $\bar{c}_m=0$, we obtain $\Breve{\Phi}_m=\sum_{j=1}^m (\tilde{c}_j + \bar{c}_j) \phi_j$, which proves our assert since $c^n_j(\triangleq \tilde{c}_j + \bar{c}_j)$ is a function of $\mathcal{L}_Gh, \mathcal{L}^2_Gh, \cdots, \mathcal{L}^m_Gh$ $\blacksquare$}

\vskip.2cm

\noindent It also holds the following result:

\begin{lm}\label{LemmaRisB}
If $w\neq 0$, $\hat{\phi}_m=\sum_{j=1}^m b^n_j(\mathcal{L}_Gh, \mathcal{L}^2_Gh, \cdots, \mathcal{L}^m_Gh) \Phi_j$, i.e., the vector $\hat{\phi}_m$ is a linear combination of the vectors  $\Phi_j$ ($j=1,\cdots,m$), where the coefficients ($b^n_j$) depend on the state only through the functions that generate the codistribution $L^m$
\end{lm}

\proof{We proceed by induction. By a direct computation it is immediate to obtain: $\hat{\phi}_1 = \Phi_1  \frac{1}{\mathcal{L}_Gh}$ (note that $\mathcal{L}_Gh=L^1_gw\neq 0$). 

\noindent {\bf Inductive step:} Let us assume that $\hat{\phi}_{m-1}=\sum_{j=1}^{m-1} b_j(\mathcal{L}_Gh, \mathcal{L}^2_Gh, \cdots, \mathcal{L}^{m-1}_Gh) \Phi_j$. We need to prove that $\hat{\phi}_m=\sum_{j=1}^m b^n_j(\mathcal{L}_Gh, \mathcal{L}^2_Gh, \cdots, \mathcal{L}^m_Gh) \Phi_j$. 
We start by applying on both members of the equality $\hat{\phi}_{m-1}=\sum_{j=1}^{m-1} b_j(\mathcal{L}_Gh, \mathcal{L}^2_Gh, \cdots, \mathcal{L}^{m-1}_Gh) \Phi_j$ the Lie bracket with respect to $G$. We obtain for the first member: $[\hat{\phi}_{m-1}, ~G]= \hat{\phi}_m \mathcal{L}^1_Gh$. For the second member we have:

\[
\sum_{j=1}^{m-1} [b_j \Phi_j, ~G]= \sum_{j=1}^{m-1} b_j [\Phi_j, ~G]- \sum_{j=1}^{m-1} \mathcal{L}_G b_j \Phi_j=
\]
\[
=\sum_{j=1}^{m-1} b_j \Phi_{j+1}- \sum_{j=1}^{m-1}  \sum_{i=1}^{m-1} \frac{\partial b_j}{\partial (\mathcal{L}_G^ih)} \mathcal{L}_G^{i+1}h\Phi_j
\]

\noindent Since $\mathcal{L}_Gh=L^1_gw\neq 0$, by setting $\tilde{b}_j=\frac{b_{j-1}}{\mathcal{L}^1_Gh}$ for $j=2,\cdots,m$ and $\tilde{b}_1=0$, and by setting $\bar{b}_j=-\sum_{i=1}^{m-1} \frac{\partial b_j}{\partial (\mathcal{L}_G^ih)} \frac{\mathcal{L}_G^{i+1}h}{\mathcal{L}^1_Gh}$ for $j=1,\cdots,m-1$ and $\bar{b}_m=0$, we obtain $\hat{\phi}_m=\sum_{j=1}^m (\tilde{b}_j + \bar{b}_j) \Phi_j$, which proves our assert since $b^n_j(\triangleq \tilde{b}_j + \bar{b}_j)$ is a function of $\mathcal{L}_Gh, \mathcal{L}^2_Gh, \cdots, \mathcal{L}^m_Gh$ $\blacksquare$}

\vskip.2cm

\noindent An important consequence of the previous two lemmas is the following result:

\begin{pr}\label{PropertyLemma2}
If $w\neq 0$, the following two codistributions coincide:
\begin{enumerate}
\item $span\{ \mathcal{L}_{\Phi_0}d  h, \mathcal{L}_{\Phi_1}d  h, \cdots, \mathcal{L}_{\Phi_m}d  h, \mathcal{L}^1_Gd  h, \cdots\mathcal{L}^m_Gd  h\}$;
\item $span\{ \mathcal{L}_{\hat{\phi}_0}d  h, \mathcal{L}_{\hat{\phi}_1}d  h, \cdots, \mathcal{L}_{\hat{\phi}_m}d  h, \mathcal{L}^1_Gd  h, \cdots\mathcal{L}^m_Gd  h\}$;
\end{enumerate}
\end{pr}

\noindent We are now ready to prove the following fundamental result:

\begin{theorem}[Separation]\label{TheoremSeparation}
If $w\neq 0$, $\bar{\Omega}_m=\tilde{\Omega}_m\triangleq [\Omega_m,0_k]+ L^m$
\end{theorem}

\proof{We proceed by induction. By definition, $\bar{\Omega}_0=\tilde{\Omega}_0$ since they are both the span of $d  h$.

\noindent {\bf Inductive step:} Let us assume that $\bar{\Omega}_{m-1}=\tilde{\Omega}_{m-1}$. We have:
$\bar{\Omega}_m=\bar{\Omega}_{m-1}+\mathcal{L}_F \bar{\Omega}_{m-1} + \mathcal{L}_G \bar{\Omega}_{m-1}=\bar{\Omega}_{m-1} + \mathcal{L}_F \tilde{\Omega}_{m-1} + \mathcal{L}_G \bar{\Omega}_{m-1}=\bar{\Omega}_{m-1} + [\mathcal{L}_f \Omega_{m-1}, 0_k]+\mathcal{L}_F L^{m-1} + \mathcal{L}_G \bar{\Omega}_{m-1}$.
On the other hand, $\mathcal{L}_F L^{m-1} = \mathcal{L}_F \mathcal{L}^1_Gd  h + \cdots + \mathcal{L}_F \mathcal{L}^{m-2}_Gd  h + \mathcal{L}_F \mathcal{L}^{m-1}_Gd  h$. The first $m-2$ terms are in $\bar{\Omega}_{m-1}$. Hence we have: $\bar{\Omega}_m= \bar{\Omega}_{m-1} + [\mathcal{L}_f \Omega_{m-1},0_k] + \mathcal{L}_F \mathcal{L}^{m-1}_G d  h + \mathcal{L}_G \bar{\Omega}_{m-1}$. By using lemma \ref{LemmaLemma1} we obtain: $\bar{\Omega}_m= \bar{\Omega}_{m-1} + [\mathcal{L}_f \Omega_{m-1},0_k] + \mathcal{L}_{\Phi_{m-1}} d  h + \mathcal{L}_G \bar{\Omega}_{m-1}$. By using again the induction assumption we obtain: $\bar{\Omega}_m= [\Omega_{m-1},0_k] + L^{m-1} + [\mathcal{L}_f \Omega_{m-1},0_k] + \mathcal{L}_{\Phi_{m-1}} d  h + \mathcal{L}_G [\Omega_{m-1},0_k] + \mathcal{L}_G L^{m-1}=[\Omega_{m-1},0_k] + L^m + [\mathcal{L}_f \Omega_{m-1},0_k] + \mathcal{L}_{\Phi_{m-1}} d  h + [\mathcal{L}_{\frac{g}{L^1_g}} \Omega_{m-1},0_k]$, where we used $L^m+\mathcal{L}_G [\Omega_{m-1},0_k]=L^m+ [\mathcal{L}_{\frac{g}{L^1_g}} \Omega_{m-1},0_k]$, which holds because $\mathcal{L}_Gh=L^1_gw\neq 0$.
By using proposition \ref{PropertyLemma2}, we obtain: $\bar{\Omega}_m=[\Omega_{m-1},0_k] + L^m + [\mathcal{L}_f \Omega_{m-1},0_k] + \mathcal{L}_{\hat{\phi}_{m-1}} d  h + [\mathcal{L}_{\frac{g}{L^1_g}} \Omega_{m-1},0_k]=\tilde{\Omega}_m$ $\blacksquare$}

\noindent Theorem \ref{TheoremSeparation} is fundamental. It allows us to obtain all the observability properties of the original state by restricting the computation to the codistribution defined by algorithm \ref{AlgoO1}, namely a codistribution whose covectors have the same dimension of the original space. In other words, the dimension of these covectors is independent of the state augmentation. 

\subsection{Convergence}\label{SubSectionConvergence}

Algorithm \ref{AlgoO1} is recursive and $\Omega_m\subseteq\Omega_{m+1}$. This means that, if for a given $m$ the differentials of the components of the original state belong to $\Omega_m$, we can conclude that the original state is weakly locally observable. On the other hand, if this is not true, we cannot exclude that it is true for a larger $m$. The goal of this section is precisely to address this issue. We will show that the algorithm converges in a finite number of steps and we will also provide the criterion to establish that the algorithm has converged (theorem \ref{TheoremStop}). This theorem will be proved at the end of this section since we need to introduce several important new quantities and properties.

\noindent When investigating the convergence properties of algorithm \ref{AlgoO1}, we remark that, the main difference between algorithm \ref{AlgoHK} and \ref{AlgoO1}, is the presence of the last term in the recursive step of the latter. Without this term, the convergence criterion would simply consist of the inspection of the equality $\Omega_{m+1}=\Omega_m$, as for algorithm \ref{AlgoHK}. Indeed, the convergence criterion of algorithm \ref{AlgoHK} is a consequence of the fact that, all the terms that appear at the recursive step of algorithm \ref{AlgoHK}, are the Lie derivative of the codistribution at the previous step, along fixed vector fields (i.e., vector fields that remain the same at each step of the algorithm). This is not the case for the last term at the recursive step of algorithm \ref{AlgoO1}.

The following result provides the convergence criterion in a very special case that basically occurs when the contribution due to the last term in the recursive step of algorithm \ref{AlgoO1} is included in the other terms. In this case, we obviously obtain that the convergence criterion consists of the inspection of the equality $\Omega_{m+1}=\Omega_m$, as for algorithm \ref{AlgoHK}.
For any integer $j\ge 0$ we define:

\begin{equation}\label{EquationDefChi}
\chi_j \triangleq \frac{\mathcal{L}_{\phi_j}L^1_g}{L^1_g}
\end{equation}

\noindent We have the following result:

\begin{lm}\label{LemmaConvergenceSpecial}
Let us denote by $\Lambda_j$ the distribution generated by $\phi_0, \phi_1, \cdots, \phi_j$ and by $m(\le n-1)$ the smallest integer for which $\Lambda_{m+1}=\Lambda_m$ ($n$ is the dimension of the state $x$). In the very special case when $\chi_j=0$, $\forall j=0,\cdots,m$, algorithm \ref{AlgoO1} converges at the integer $j$ such that $\Omega_{j+1}=\Omega_j$ and this occurs in at most $n-1$ steps\footnote{Note that, in the case of multiple known inputs (i.e., $m_u>1$), the previous result does not change. In particular, $\Lambda_j$ becomes the distribution generated by $^1\phi_0, ^1\phi_1, \cdots, ^1\phi_j, \cdots, ^{m_u}\phi_0, ^{m_u}\phi_1, \cdots, ^{m_u}\phi_j$ and the very special case is now characterized by 
$\frac{\mathcal{L}_{^i\phi_j}L^1_g}{L^1_g}=0$, $\forall j=0,\cdots,m$ and $\forall i=1,\cdots,m_u$.}.
\end{lm}

\proof{First of all, we remind the reader that the existence of an integer $m(\le n-1)$ such that $\Lambda_{m+1}=\Lambda_m$ is proved in \cite{Isi95}. In particular, the first chapter in \cite{Isi95} analyzes the convergence of $\Lambda_j$ with respect to $j$. It is proved that the distribution converges to $\Lambda^*$ and that the convergence is achieved  at the smallest integer for which we have $\Lambda_{m+1}=\Lambda_m$. Additionally, we have $\Lambda_{m+1}=\Lambda_m=\Lambda^*$ and $m$ cannot exceed $n-1$.

In the very special case when $\chi_j=0$, $\forall j=0,\cdots,m$, thanks to the aforementioned convergence of the distribution $\Lambda_j$, we easily obtain that $\mathcal{L}_{\phi_{j-1}}\mathcal{L}_gh=0$ $\forall j\ge 1$.
Now, let us consider the following equation:

\begin{equation}\label{EquationLphi_jh_1}
\mathcal{L}_{\phi_j}h=\frac{1}{L^1_g}\left( \mathcal{L}_{\phi_{j-1}}\mathcal{L}_g h-\mathcal{L}_g\mathcal{L}_{\phi_{j-1}}h\right)
\end{equation}

\noindent Since $\mathcal{L}_{\phi_{j-1}}\mathcal{L}_gh=0$ $\forall j\ge 1$, we have $\mathcal{L}_{\phi_j}h=-\mathcal{L}_{\frac{g}{L^1_g}}\mathcal{L}_{\phi_{j-1}}h$, for any $j\ge 1$. Therefore, we conclude that, the last term in the recursive step of algorithm \ref{AlgoO1}, is included in the second last term and,  in this special case, algorithm \ref{AlgoO1} has converged when $\Omega_{m+1}=\Omega_m$. This occurs in at most $n-1$ steps, as for algorithm \ref{AlgoHK}. 
$\blacksquare$}

\vskip .2cm
\noindent Let us consider now the general case. To proceed we need to introduce several important new quantities and properties.

\noindent For a given positive integer $j$ we define the vector $\psi_j\in \mathbb{R}^n$ by the following algorithm:

\begin{enumerate}
\item $\psi_0=f$;
\item $\psi_j=[\psi_{j-1}, ~\frac{g}{L^1_g}]$
\end{enumerate}

\noindent It is possible to find the expression that relates these vectors to the vectors $\phi_j$, previously defined. Specifically we have:

\begin{lm}\label{LemmaPsiPhi}
It holds the following equation:

\begin{equation}\label{EquationPsiPhi}
\psi_j=\phi_j + \left\{ \sum_{i=0}^{j-1} (-)^{j-i} \mathcal{L}^{j-i-1}_{\frac{g}{L^1_g}} \left( \frac{\mathcal{L}_{\phi_i}L^1_g}{L^1_g}\right) \right\} \frac{g}{L^1_g}
\end{equation}

\end{lm}

\proof{We proceed by induction. By definition $\psi_0=\phi_0=f$ and equation (\ref{EquationPsiPhi}) holds for $j=0$.  

\noindent {\bf Inductive step:} Let us assume that it holds for a given $j-1\ge 0$ and let us prove its validity for $j$. We have:

\[
\psi_j=\left[\psi_{j-1}, ~\frac{g}{L^1_g}\right]=\left[\phi_{j-1}, ~\frac{g}{L^1_g}\right]\]
\[+ \left[\left\{ \sum_{i=0}^{j-2} (-)^{j-i-1} \mathcal{L}^{j-i-2}_{\frac{g}{L^1_g}} \left( \frac{\mathcal{L}_{\phi_i}L^1_g}{L^1_g}\right) \right\} \frac{g}{L^1_g}, ~\frac{g}{L^1_g}\right]
\]

\noindent On the other hand:

\[
\left[\phi_{j-1}, ~\frac{g}{L^1_g}\right] = \phi_{j} - \frac{\mathcal{L}_{\phi_{j-1}}L^1_g}{L^1_g} \frac{g}{L^1_g}
\]

\noindent and

\[
\left[\left\{ \sum_{i=0}^{j-2} (-)^{j-i-1} \mathcal{L}^{j-i-2}_{\frac{g}{L^1_g}} \left( \frac{\mathcal{L}_{\phi_i}L^1_g}{L^1_g}\right) \right\} \frac{g}{L^1_g}, ~\frac{g}{L^1_g}\right]=\]
\[-\mathcal{L}_{\frac{g}{L^1_g}}
\left\{ \sum_{i=0}^{j-2} (-)^{j-i-1} \mathcal{L}^{j-i-2}_{\frac{g}{L^1_g}} \left( \frac{\mathcal{L}_{\phi_i}L^1_g}{L^1_g}\right) \right\} \frac{g}{L^1_g}=
\]
\[
\left\{ \sum_{i=0}^{j-2} (-)^{j-i} \mathcal{L}^{j-i-1}_{\frac{g}{L^1_g}} \left( \frac{\mathcal{L}_{\phi_i}L^1_g}{L^1_g}\right) \right\} \frac{g}{L^1_g}
\]

\noindent Hence:

\[
\psi_j= \phi_{j} - \frac{\mathcal{L}_{\phi_{j-1}}L^1_g}{L^1_g} \frac{g}{L^1_g} +
\left\{ \sum_{i=0}^{j-2} (-)^{j-i} \mathcal{L}^{j-i-1}_{\frac{g}{L^1_g}} \left( \frac{\mathcal{L}_{\phi_i}L^1_g}{L^1_g}\right) \right\} \frac{g}{L^1_g},
\]

\noindent which coincides with (\ref{EquationPsiPhi})
$\blacksquare$}

\vskip .2cm

\noindent We have the following result:

\begin{lm}\label{LemmaFunctionsInOmegam}
For $i=0,1,\cdots,m-2$, we have:
\begin{equation}\label{EquationFunctionsInOmegam}
\mathcal{D} \frac{\mathcal{L}_{\phi_i}L^1_g}{L^1_g} \in \Omega_m
\end{equation}

\end{lm}

\proof{By construction, $\mathcal{D} \mathcal{L}_{\phi_i}h \in \Omega_m$, for any $i=1,\cdots,m-1$. On the other hand, we have:

\[
 \mathcal{L}_{\phi_i}h=\frac{1}{L^1_g}[ \mathcal{L}_{\phi_{i-1}} \mathcal{L}_gh - \mathcal{L}_g\mathcal{L}_{\phi_{i-1}}h]= \frac{\mathcal{L}_{\phi_{i-1}}L^1_g}{L^1_g}- \mathcal{L}_{\frac{g}{L^1_g}}\mathcal{L}_{\phi_{i-1}}h
\]

We compute the differential of both members of this equation. Since $\mathcal{D} \mathcal{L}_{\frac{g}{L^1_g}}\mathcal{L}_{\phi_{i-1}}h \in \Omega_m$, for any $i=1,\cdots,m-1$, also $\mathcal{D} \frac{\mathcal{L}_{\phi_{i-1}}L^1_g}{L^1_g} \in \Omega_m$
 $\blacksquare$}

\noindent From lemma \ref{LemmaPsiPhi} with $j=1,\cdots,m-1$ and lemma \ref{LemmaFunctionsInOmegam} it is immediate to obtain the following result:

\begin{pr}\label{PropPhimTot}
If $\Omega_m$ is invariant with respect to $\mathcal{L}_f$ and $\mathcal{L}_{\frac{g}{L^1_g}}$ then it is also invariant with respect to $\mathcal{L}_{\phi_j}$, $j=1,\cdots,m-1$. 
\end{pr}


\vskip .5cm

\noindent  In order to obtain the convergence criterion for algorithm \ref{AlgoO1} we need to find the analytic expression that relates three consecutive Lie derivatives, $\mathcal{L}_{\phi_j}h$, $\mathcal{L}_{\phi_{j-1}}h$ and $\mathcal{L}_{\phi_{j-2}}h$. This will allow us to detect the key quantity that governs the convergence of algorithm \ref{AlgoO1}, in particular regarding the contribution due to the last term in the recursive step. This quantity is a scalar and it is the one provided in (\ref{EquationRho}). For the sake of clarity, we provide equation (\ref{EquationRho}) below:

\[
\tau \triangleq \frac{L^2_g}{(L^1_g)^2}
\]

\noindent where $L^2_g \triangleq \mathcal{L}^2_g h$. The behaviour of the last term in the recursive step of algorithm \ref{AlgoO1} is given by the following lemma:

\begin{lm}\label{LemmaKeyEquality}
We have the following key equality:

\begin{equation}\label{EquationKeyEquality}
\mathcal{L}_{\phi_j}h=\mathcal{L}_{\phi_{j-2}} \tau + \tau \frac{\mathcal{L}_{\phi_{j-2}}L^1_g}{L^1_g}-\mathcal{L}_{\frac{g}{L^1_g}}\left(\frac{\mathcal{L}_{\phi_{j-2}}L^1_g}{L^1_g}+\mathcal{L}_{\phi_{j-1}}h \right)
\end{equation}
\noindent $j\ge 2$.
\end{lm}

\proof{We will prove this equality by an explicit computation. We have:

\[
\mathcal{L}_{\phi_j}h=\frac{1}{L^1_g}\left( \mathcal{L}_{\phi_{j-1}}\mathcal{L}_g h-\mathcal{L}_g\mathcal{L}_{\phi_{j-1}}h\right)
\]

The second term on the right hand side simplifies with the last term in (\ref{EquationKeyEquality}). Hence we have to prove:

\begin{equation}\label{EquationKeyEqualityProof}
\frac{1}{L^1_g}\mathcal{L}_{\phi_{j-1}}L^1_g=\mathcal{L}_{\phi_{j-2}} \tau + \tau \frac{\mathcal{L}_{\phi_{j-2}}L^1_g}{L^1_g}-\mathcal{L}_{\frac{g}{L^1_g}}\frac{\mathcal{L}_{\phi_{j-2}}L^1_g}{L^1_g}
\end{equation}

We have:

\begin{equation}\label{EquationKeyEqualityProof1}
\frac{1}{L^1_g}\mathcal{L}_{\phi_{j-1}}L^1_g=\frac{1}{(L^1_g)^2}\left( \mathcal{L}_{\phi_{j-2}} L^2_g-\mathcal{L}_g\mathcal{L}_{\phi_{j-2}}L^1_g\right)
\end{equation}

We remark that:

\[
\frac{1}{(L^1_g)^2} \mathcal{L}_{\phi_{j-2}} L^2_g=\mathcal{L}_{\phi_{j-2}} \tau + 2 \tau \frac{\mathcal{L}_{\phi_{j-2}}L^1_g}{L^1_g}
\]

and

\[
\frac{1}{(L^1_g)^2}\mathcal{L}_g\mathcal{L}_{\phi_{j-2}}L^1_g = \tau \frac{\mathcal{L}_{\phi_{j-2}}L^1_g}{L^1_g} + \mathcal{L}_{\frac{g}{L^1_g}}\frac{\mathcal{L}_{\phi_{j-2}}L^1_g}{L^1_g}
\]

By substituting these two last equalities in (\ref{EquationKeyEqualityProof1}) we immediately obtain (\ref{EquationKeyEqualityProof}) $\blacksquare$}

\begin{lm}\label{LemmaRhoInOmegam}
In general, it exists an integer $m\le n+2$ (being $n$ the dimension of $x$) such that $\mathcal{D}\tau \in \Omega_m$.
\end{lm}

\proof{Let us introduce the following notation, for a given integer $j$:

\begin{itemize}

\item $\mathcal{Z}_j\triangleq  \mathcal{L}_{\phi_{j+2}}h$;

\item $\mathcal{B}_j \triangleq \mathcal{L}_{\phi_j}\tau$;

\item $\chi_j \triangleq \frac{\mathcal{L}_{\phi_j}L^1_g}{L^1_g}$.

\end{itemize}
 
\noindent By construction, $\mathcal{D} \mathcal{Z}_j \in \Omega_{j+3}$. On the other hand, from equation (\ref{EquationKeyEquality}), we immediately obtain: 

\begin{equation}\label{EquationZizza}
\mathcal{D} \mathcal{Z}_j = \mathcal{D}  \mathcal{B}_j + \chi_j \mathcal{D}\tau + \tau \mathcal{D}\chi_j-
\mathcal{L}_{\frac{g}{L^1_g}}\left(\mathcal{D}\chi_j + \mathcal{D} \mathcal{L}_{\phi_{j+1}}h\right)
\end{equation}

\noindent By using lemma \ref{LemmaFunctionsInOmegam} we obtain the following results:

\begin{itemize}

\item $\tau \mathcal{D}\chi_j \in \Omega_{j+2}$;

\item $\mathcal{L}_{\frac{g}{L^1_g}}\mathcal{D}\chi_j \in \Omega_{j+3}$.

\end{itemize}

\noindent Additionally, $\mathcal{L}_{\frac{g}{L^1_g}}\mathcal{D} \mathcal{L}_{\phi_{j+1}}h  \in \Omega_{j+3}$. Hence, from (\ref{EquationZizza}), we obtain that the following covector: 

\begin{equation}\label{EquationZizzap}
\mathcal{Z}_j' \triangleq \mathcal{D}\mathcal{B}_j +\chi_j \mathcal{D}\tau
\end{equation}

\noindent belongs to $ \Omega_{j+3}$. Let us denote by $j^*$ the smallest integer such that:

\begin{equation}\label{Equationj*}
\mathcal{D}\mathcal{B}_{j^*}=\sum_{j=0}^{j^*-1} c_j \mathcal{D}\mathcal{B}_j + c_{-1}\mathcal{D} h
\end{equation}

\noindent Note that $j^*$ is a finite integer and in particular $j^*\le n-1$. Indeed, if this would not be the case, the dimension of the codistribution generated by $\mathcal{D}h, \mathcal{D}\mathcal{B}_0,\mathcal{D}\mathcal{B}_1, \cdots, \mathcal{D}\mathcal{B}_ {n-1}$ would be $n+1$, i.e., larger than $n$.
From (\ref{Equationj*}) and (\ref{EquationZizzap}) we obtain:

\begin{equation}\label{EquationZizza*}
\mathcal{Z}_{j*}' = \sum_{j=0}^{j^*-1} c_j \mathcal{D}\mathcal{B}_j + c_{-1}\mathcal{D} h+\chi_{j*} \mathcal{D} \tau
\end{equation}

\noindent From equation (\ref{EquationZizzap}), for $j=0,\cdots,j^*-1$, we obtain: $\mathcal{D}\mathcal{B}_j= \mathcal{Z}_j' - \chi_j \mathcal{D} \tau$.  By substituting in (\ref{EquationZizza*}) we obtain: 

\begin{equation}
\mathcal{Z}_{j*}'-\sum_{j=0}^{j^*-1} c_j \mathcal{Z}_j'  - c_{-1}\mathcal{D} h = \left(-\sum_{j=0}^{j^*-1} c_j \chi_j + \chi_{j*}  \right) \mathcal{D} \tau
\end{equation}

\noindent We remark that the left hand side consists of the sum of covectors that belong to $\Omega_{j^*+3}$. Since in general $\chi_{j*}\neq \sum_{j=0}^{j^*-1} c_j\chi_j$, we have $\mathcal{D}\tau \in \Omega_{j^*+3}$. By setting $m\triangleq j^*+3$, we have $m\le n+2$ and $\mathcal{D}\tau \in \Omega_m$
$\blacksquare$}

\noindent The previous lemma ensures that, in general, it exists a finite $m\le n+2$ such that $\mathcal{D}\tau \in \Omega_m$. Note that the previous proof holds if the quantity $\chi_{j*}- \sum_{j=0}^{j^*-1} c_j\chi_j$ does not vanish. This holds in general, with the exception of the trivial case considered in lemma \ref{LemmaConvergenceSpecial}, in which case $\chi_j=0, ~\forall j$.

 The following theorem allows us to obtain the criterion to stop algorithm \ref{AlgoO1}:


\begin{theorem}\label{TheoremStop}
If $\mathcal{D}\tau \in \Omega_m$ and $\Omega_m$ is invariant under $\mathcal{L}_f$ and $\mathcal{L}_{\frac{g}{L^1_g}}$, then $\Omega_{m+p}=\Omega_m$ $\forall p\ge 0$
\end{theorem}

\proof{We proceed by induction. Obviously, the equality holds for $p=0$. 

\noindent {\bf Inductive step:} let us assume that $\Omega_{m+p}=\Omega_m$ and let us prove that $\Omega_{m+p+1}=\Omega_m$. We have to prove that $\mathcal{D}\mathcal{L}_{\phi_{m+p}}h \in \Omega_m$. Indeed, from the inductive assumption, we know that $\Omega_{m+p}(=\Omega_m)$ is invariant under $\mathcal{L}_f$ and $\mathcal{L}_{\frac{g}{L^1_g}}$. Additionally, because of this invariance, by using proposition \ref{PropPhimTot}, we obtain that $\Omega_m$ is also invariant under $\mathcal{L}_{\phi_j}$, for $j=1,2,\cdots,m+p-1$. Since $\mathcal{D}\tau \in \Omega_m$ we have $\mathcal{D}\mathcal{L}_{\phi_{m+p-2}} \tau \in \Omega_m$. Additionally,  $\mathcal{D}\mathcal{L}_{\phi_{m+p-1}}h \in \Omega_m$ and, because of lemma \ref{LemmaFunctionsInOmegam}, we also have $\mathcal{D}\frac{\mathcal{L}_{\phi_{m+p-2}}L^1_g}{L^1_g} \in \Omega_m$. Finally, because of the invariance under $\mathcal{L}_{\frac{g}{L^1_g}}$, also the Lie derivatives along $\frac{g}{L^1_g}$ of $\mathcal{D}\mathcal{L}_{\phi_{m+p-1}}h$ and $\mathcal{D}\frac{\mathcal{L}_{\phi_{m+p-2}}L^1_g}{L^1_g}$ belong to $\Omega_m$.
Now, we use equation (\ref{EquationKeyEquality}) for $j=m+p$. By computing the differential of this equation it is immediate to obtain that $\mathcal{D}\mathcal{L}_{\phi_{m+p}}h \in \Omega_m$
$\blacksquare$}

\vskip .3cm

We conclude this section by providing an upper bound for the number of steps that are in general necessary to achieve the convergence. The dimension of $\Omega_{j^*+2}$ is at least the dimension of the span of the covectors: $\mathcal{D}h, ~ \mathcal{Z}'_0, ~ \mathcal{Z}'_1, ~\cdots, ~ \mathcal{Z}'_{j^*-1}$. From the definition of $j^*$, we know that the vectors $\mathcal{D}h, ~\mathcal{D} \mathcal{B}_0, ~\mathcal{D} \mathcal{B}_1, ~\cdots, ~\mathcal{D} \mathcal{B}_{j^*-1}$ are independent meaning that the dimension of their span is $j^*+1$. Hence, from (\ref{EquationZizzap}), it easily follows that the dimension of the span of the vectors  $\mathcal{D}h, ~ \mathcal{Z}'_0, ~ \mathcal{Z}'_1, ~\cdots, ~ \mathcal{Z}'_{j^*-1}, ~\mathcal{D}\tau$ is at least $j^*+1$. Since $\Omega_{j^*+3}$ contains this span, its dimension is at least $j^*+1$. Therefore, the condition $\Omega_{m+1}=\Omega_m$, for $m\ge j^*+3$ is achieved for $m\le n+2$.

\subsection{Extension to the case of multiple known inputs}\label{SubSectionExt}

It is immediate to repeat all the steps carried out in the previous two subsections and extend the validity of theorem \ref{TheoremSeparation} to the case of multiple known inputs ($m_u>1$).
Additionally, also theorem \ref{TheoremStop} can be easily extended to cope with the case of multiple known inputs. In this case, requiring that $\Omega_{m+1} =\Omega_m$ means that $\Omega_m$ must be invariant with respect to $\mathcal{L}_{\frac{g}{L^1_g}}$ and all $\mathcal{L}_{f^i}$ simultaneously.

\chapter{Conclusion}\label{ChapterConclusion}

The goal of this work was to present the analytic solution of a fundamental open problem in control theory.
The problem is known in the literature as the Unknown Input Observability (UIO) problem. It consists in deriving the analytic criterion that allows us to check the system observability in presence of unknown inputs (throughout this work we also used the term {\it disturbance} to mean an unknown input).
The solution here provided holds in the case of a single unknown input.
In other words, it is provided the analytic criterion that allows us to obtain the observability of a nonlinear system in presence of a single unknown input and multiple known inputs (see chapter \ref{ChapterSystemsDefinition} for a definition of the considered systems).

As for the observability rank condition, the proposed analytic criterion is based on the computation of the observable codistribution. Similarly to the case of only known inputs, the observable codistribution is obtained by recursively computing the Lie derivatives of the outputs along the vector fields that characterize the dynamics. However, in correspondence of the unknown input, the corresponding vector field must be suitably rescaled. Additionally, the Lie derivatives of the outputs must also be computed along a new set of vector fields that are obtained by recursively performing suitable Lie bracketing of the vector fields that define the dynamics. In practice, the entire observable codistribution is obtained by a very simple recursive algorithm. This algorithm was provided in chapter \ref{ChapterObsCod}. Finally, it was shown that the recursive algorithm converges in a finite number of steps and the criterion to establish that the convergence has been reached was provided.

The analytic criterion here introduced  is the extension of the observability rank condition to the case when the dynamics are driven by an unknown input (and multiple known inputs). The analytic criterion was quickly outlined in chapter \ref{ChapterEORC} to make easy its automatic implementation. Its validity was proven in chapter \ref{ChapterProofs}.
Finally, it  was illustrated in chapter \ref{ChapterApplications} by checking the weak local observability of several nonlinear systems driven by multiple known inputs and a single unknown input.

We are currently extending the proposed analytic criterion to deal with the general case of multiple unknown inputs and when the dynamics are affine (and not simply linear) in the inputs (both known and unknown). In few words, to deal with the general case described by equation (\ref{EquationSystemDefinitionDynamicsG}).

We conclude this work with the following remark, which has a fundamental practical importance for designing an estimator in presence of a single unknown input. We remark that, the state estimated by the two estimators adopted in chapter \ref{ChapterApplications} (sections \ref{SubSectionSimulations2D} and \ref{SubSectionSimulations3D}), includes the disturbance together with its first order time derivative. We found a much worse performance by using an estimator that only estimates the original state or by including in the state only the disturbance (without its time derivative). 
%
%
%
%
%
%
%

We conjecture that this is a consequence of the fact that the observable codistibution is $\Omega_2$, i.e., for the considered system, algorithm \ref{AlgoO1} converges at its second step. 
Indeed, we proved that, for a system for which algorithm \ref{AlgoO1} converges to $\Omega_m$ (i.e., converges in $m$ steps) the observable modes\footnote{See definition \ref{DefObservableMode} for the definition of {\it observable mode}.} satisfy a fundamental property. These functions are constant on the indistinguishable sets only in the extended state space that includes the disturbance and its time derivatives up to the $m-1$ order (see proposition \ref{PrConstantLieDer} together with theorem \ref{TheoremSeparation} and lemma \ref{LemmaFundamental}). Note that, in the space of the original state (or in a less extended space, namely, that does not include all the first $m-1$ time derivatives), the observable modes are not constant on the indistinguishable sets. Proving this conjecture will require further theoretical investigations. The validity of this conjecture would make the analytic tool not only a tool to check the state observability but also a fundamental tool to design an appropriate estimator for a given nonlinear system in presence of a disturbance.

\appendix

\chapter{The case when $L^1_g=0$}\label{AppendixCanonic}

We consider the case of multiple outputs. In other words, our system is defined by the following equations:

\begin{equation}\label{EquationAppendixCanonicG0}
\left\{\begin{aligned}
  \dot{x} &=  \sum_{i=1}^{m_u}f^i ( x ) u_i + g ( x ) w \\
  y_k &=h_k(x),  ~~~ k=1,\cdots, p\\
\end{aligned}\right.
\end{equation}

\noindent We introduce the following notation:

\begin{itemize}

\item $\mathcal{F}$ is the space of functions that contains all the outputs $h_1,\cdots,h_p$ and their Lie derivatives along the vector fields $f^1,\cdots,f^{m_u}$, up to any order;

\item $\mathcal{D} \mathcal{F}$ is the codistribution generated by the differentials, with respect to the state $x$, of the functions in $\mathcal{F}$;

\item $\mathcal{L}_G \mathcal{F}$ is the space of functions that consists of the Lie derivatives along $G$ of the functions in $\mathcal{F}$ (we remind the reader that $G$ is the vector defined in (\ref{EquationFG}));


\end{itemize}

\noindent It is immediate to prove the following properties:

\begin{itemize}

\item From $\mathcal{F}$ we can select a set of functions such that their differentials generate $\mathcal{D} \mathcal{F}$.

\item It exists an integer $m$ such that, by running $m$ recursive steps of algorithm \ref{AlgoOT}, we obtain a codistribution $\bar{\Omega}_m$ that contains $[\mathcal{D} \mathcal{F},~0_{m}]$ (i.e., the codistribution $\mathcal{D} \mathcal{F}$ once embedded in $\mathds{R}^{n+m}$)\footnote{Note that, from proposition \ref{Pr3.1}, the functions in $\mathcal{F}$ are weakly observable.}.

\item In general, the functions belonging to $\mathcal{L}_G  \mathcal{F}$, are functions of $x$ and $w$.


\end{itemize}

\noindent Because of the second property, we are allowed to choose any function in $\mathcal{F}$ as a system output (note that the functions in $\mathcal{F}$ only depend on the state $x$ and not on the unknown input $w$).
Hence, if for a given system the condition $\mathcal{L}_g h_k=0$, $\forall k=1,\cdots p$, we can redefine the output by choosing any function in $\mathcal{F}$.

To conclude, let us consider the case when no function in $\mathcal{L}_G  \mathcal{F}$ depends on $w$.
This means that $\mathcal{L}_G \mathcal{F}$ only contains the zero function (the Lie derivative along $G$ of any function in $\mathcal{F}$ vanishes). As a result, any order Lie derivative of any outputs, computed along $F^1,\cdots,F^{m_u}$ (which are the vector fields defined in (\ref{EquationFG})) and at least once along $G$, vanishes. This means that, in algorithm \ref{AlgoOT}, we can ignore the contribution due to the Lie derivative along $G$. But this trivially means that the observable codistribution is precisely $\mathcal{D} \mathcal{F}$ and the unknown input is spurious.

\end{document}